\documentclass[12pt,letterpaper]{article}
\usepackage{amssymb}
\usepackage{amsmath}
\usepackage{amscd}
\usepackage{graphicx}
\usepackage{tikz}
\usepackage{lscape}
\usepackage{arydshln}

\def\be{\begin{equation}}
\def\ee{\end{equation}}
\def\bea{\begin{eqnarray}}
\def\eea{\end{eqnarray}}
\def\bse{\begin{subequations}}
\def\ese{\end{subequations}}

\def\vCent#1{\vcenter{\hbox{\hss#1\hss}}}

\def\bsf{\sffamily\bfseries}

\def\4{\text{\bsf-}}

\def\hi{{\hat\imath}}

\def\hj{{\hat\jmath}}

\def\fc#1#2{\relax\ifmmode{\scriptstyle\frac{#1}{#2}} 
                    \else$\scriptstyle\frac{#1}{#2}$\fi}    

\definecolor{Hey}{rgb}{.9,.05,.4}
\definecolor{orange}{rgb}{1,.5,0}
\definecolor{plum}{rgb}{.4,0,.6}
\definecolor{R}{rgb}{1,0,0}
\definecolor{G}{rgb}{0,1,0}
\definecolor{B}{rgb}{0,0,1}
\definecolor{Red}    {rgb}{1.00,0.00,0.00} 
\definecolor{Green}  {rgb}{0.00,0.75,0.00} 
\definecolor{Blue}   {rgb}{0.00,0.00,1.00} 
\definecolor{Orange} {rgb}{1.00,0.67,0.00} 
\definecolor{Purple} {rgb}{0.50,0.00,0.50} 
\definecolor{Gold}   {rgb}{1.00,0.90,0.00} 
\definecolor{Magenta}{rgb}{1.00,0.00,1.00} 
\definecolor{Turque} {rgb}{0.00,0.90,0.90} 
\definecolor{Seaweed}{rgb}{0.00,0.25,0.00} 
\definecolor{Brown}  {rgb}{0.50,0.13,0.00} 
\definecolor{Cobalt} {rgb}{0.00,0.00,0.50} 
\definecolor{Sage}   {rgb}{0.00,0.50,0.38} 
\definecolor{grey1}  {rgb}{0.20,0.20,0.20} 
\definecolor{grey2}  {rgb}{0.40,0.40,0.40} 
\definecolor{grey3}  {rgb}{0.60,0.60,0.60} 
\definecolor{grey4}  {rgb}{0.80,0.80,0.80} 
\definecolor{grey5}  {rgb}{0.90,0.90,0.90} 
\def\C#1#2{{\ifcase#1\or
             \color{Red}\or\color{Green}\or\color{Blue}\or
              \color{Orange}\or\color{Purple}\or\color{Gold}\or
             \color{Magenta}\or\color{Turque}\or\color{Seaweed}\or
               \color{Brown}\or\color{Cobalt}\or\color{Sage}\or
                 \color{grey1}\or\color{grey2}\or\color{grey3}\or
                 \color{grey4}\else\color{grey5}\fi#2}}
\definecolor{gray}{rgb}{.7,.7,.7}
\def\XXX{\colorbox{yellow}{\color{red}\bf X\kern-4pt{\Large$\bs*$}\kern-4.125ptX}}
\def\[{\left[}
\def\]{\right]}

\font\ro=cmsy10                          
\def\kcr{{\hbox{\ro \char'170}}}                
\def\ktl{{\hbox{\ro \char'170}}}        
\def\ktr{{\hbox{\ro \char'170}}}        
\def\kbl{{\hbox{\ro \char'170}}}        
\def\kbr{{\hbox{\ro \char'170}}}        

\parskip=\medskipamount                 
\thicklines                         

\newskip\humongous \humongous=0pt plus 1000pt minus 1000pt

\newif\ifdtup

\def\border{                                            
        \setlength{\unitlength}{1mm}
        \newcount\xco
        \newcount\yco
        \xco=-21
        \yco=12
        \begin{picture}(140,0)
        \put(\xco,\yco){$\ktl$}
        \advance\yco by-1
        {\loop
        \put(\xco,\yco){$\kcr$}
        \advance\yco by-2
        \ifnum\yco>-240
        \repeat
        \put(\xco,\yco){$\kbl$}}
        \xco=158
        \yco=12
        \put(\xco,\yco){$\ktr$}
        \advance\yco by-1
        {\loop
        \put(\xco,\yco){$\kcr$}
        \advance\yco by-2
        \ifnum\yco>-240
        \repeat
        \put(\xco,\yco){$\kbr$}}
        \put(-20,13){\tiny **University of Maryland * Center for String and
         Particle  Theory* Physics Department***University of Maryland *Center
        for String and Particle  Theory** }
        \put(-20,-241.5){\tiny **University of Maryland * Center for String and
         Particle  Theory* Physics Department***University of Maryland *Center
        for String and Particle  Theory** }
        \end{picture}
        \par\vskip-8mm}

\def\headpic{                                           
        \indent
        \setlength{\unitlength}{.4mm}
        \thinlines
        \par
        \begin{picture}(29,16)
        \put(165,16){\line(1,0){4}}
        \put(170,16){\line(1,0){4}}
        \put(180,16){\line(1,0){4}}
        \put(175,0){\line(1,0){4}}
        \put(180,0){\line(1,0){4}}
        \put(185,0){\line(1,0){4}}
        \put(169,0){\line(0,1){16}}
        \put(170,0){\line(0,1){16}}
        \put(179,0){\line(0,1){16}}
        \put(180,0){\line(0,1){16}}
        \put(184,0){\line(0,1){16}}
        \put(185,0){\line(0,1){16}}
        \put(169,16){\oval(8,32)[bl]}
        \put(170,16){\oval(8,32)[br]}
        \put(179,0){\oval(8,32)[tl]}
        \put(185,0){\oval(8,32)[tr]}
        \end{picture}
        \par\vskip-6.5mm
        \thicklines}
\def\endtitle{\end{quotation}\newpage}                  

\def\rI{{{}_{\rm I}}}
\def\rJ{{{}_{\rm J}}}

\def\hj{{\hat\jmath}}
\def\hk{{\hat k}}
\def\hi{{\hat\imath}}

\begin{document}

\border\headpic {\hbox to\hsize{Dec 2011 \hfill
{UMDEPP 011-020}
}}
\par
{$~$ \hfill
{hep-th/1201.0307}
}
\par

\setlength{\oddsidemargin}{0.3in}
\setlength{\evensidemargin}{-0.3in}
\begin{center}
\vglue .10in
{\large\bf  \center SUSY Equation Topology,\\  \vskip.05in
Zonohedra, and the Search for \\  \vskip.09in
Alternate Off-Shell Adinkras}
\\[.5in]
Keith Burghardt\footnote{keith@umd.edu}, and
S.\, James Gates, Jr.\footnote{gatess@wam.umd.edu}
\\[0.2in]

{\it Center for String and Particle Theory\\
Department of Physics, University of Maryland\\
College Park, MD 20742-4111 USA}\\[1.8in]

{\bf ABSTRACT}\\[.01in]
\end{center}
\begin{quotation}
{Results are given from a search to form adinkra-like equations
based on topologies that are not hypercubes.  An alternate class of
zonohedra topologies are used to construct adinkra-like graphs.  In particular, the rhombic dodecahedron and rhombic icosahedron are studied in detail.  Using these topological skeletons, 
equations similar to those of a supersymmetric system are found.  But 
these fail to have the interpretation of an off-shell supersymmetric system 
of equations.}

${~~~}$ \newline 
\endtitle


\section{Introduction}
$~~~~$ From the introduction of the notion of graphical representations of supersymme- trical (SUSY) theory $\cite{Adinkra}$, work on these have mostly been based on a particular class of polytopes, hypercubes, as pointed out in later work of the DFGHILM
\cite{Cubes} collaboration.  It is our ``empirical'' observation that all adinkras that
describe off-shell theories so far fall into three categories:  hypercubes,  folded 
hypercubes (as in \cite{Cubes}), or gnomons\footnote{In geometry a gnomon is defined 
as the part of a parallelogram left when a similar parallelo- \newline $~~~~\,~$ gram 
has been taken from its corner.  An analogous definition may be applied to adinkras
 \newline $~~~~\,~$ or linear combinations thereof.} 
formed from hypercubes (as in \cite{N2}).  The adinkras associated with on-shell
theories are not {\em {generally}} related to hypercubes (see \cite{Genomics}).

In order for these systems to be supersymmetric, however, any two 
super-differential operators ${\rm D}{}_{{}_{\rm I}}$ and ${\rm D}{}_{{}_{\rm J}}$, must 
obey the following property (supersymmetry algebra):
\be\label{eq:EQ1}
\{  \, {\rm D}{}_{{}_{\rm I}} ~,~ {\rm D}{}_{{}_{\rm J}}  \, \} ~=~\dot{\imath}  \, 2 \, \delta_{{}_{\rm {I \, J}}} 
 \,
 \partial_{\tau} ~~,
\ee
{\em {without}} imposing any differential equation with respect to $\tau$ on the functions in
a supermultiplet.

$~~~$The edges of any face in a polytope related to an adinkra are associated with one of 
the ${\rm D}{}_{{}_{}}$-operators in a SUSY equation.  Clearly, if we attempted to create an adinkra from a polyhedron 
where a face is described by more than two linearly independent vectors, then above could not hold.  This limits our notions of plausible adinkras to objects with \emph{only} four sided faces, i.\ e.\ faces that are completely described by {\rm {only}} two linearly independent vectors.  Further restrictions, explained below, filter out most polyhedra from even becoming adinkra candidates. Our motivation for this paper is to answer a simple question: can alternate zonohedra replace cubes in the
construction of adinkras that lead to off-shell SUSY systems? \cite{ERIK}. 
\section{Zonohedra and Alternate Polyhedra}

$~~~~$ Zonohedra are polyhedra that are symmetric with respect to 180 degree rotations\footnote{An accessible introduction to the subject of zonohedra can be found on-line via \\ $~~~~~~$
 http://en.wikipedia.org/wiki/Zonohedron.}. If we tried making adinkras from a polyhedron without this property, then two or more points could be attached to one another while both are simultaneously on the top row, or the bottom row of an associated valise adinkra graph. This is equivalent to a fermion (boson) having a fermion (boson) SUSY partner, which is invalid for an adinkra. Therefore, any polyhedron used to construct off-shell adinkras must be a zonohedron, or a higher dimensional analogue.  Higher dimensional zonotopes in  $\bf{R}^d$ with ${\bf {d >}}$ $\bf 
3$ are possible as starting points for SUSY theories with more than three ${\rm D}{}_{{}_{}}$-operators.

Not all zonohedra are qualified to be adinkras, however, as many have non-four-sided faces. Under this additional restriction for zonohedra, we get a handful of important SUSY adinkra candidates.\footnote{There is one subtlety in applying the argument above.  The actual shape of adinkras are not \\ $~~~~~$
important to the construction of the corresponding sets of equations.  Instead, it is the topol-  \\ $~~~~~$ ogy
of the geometrical object taken as the starting point in the construction of the equations
 \\ $~~~~~$
that matters.  
In other words, the relationship by which the vertices are connected one to each
 \\ $~~~~~$
other plays the critical role.}

\section{Case Study: Rhombic Dodecahedron}

$~~~$ $~~~$ Perhaps the simplest example of an adinkra candidate is the Rhombic Dodecahedron or RD (Figure$~\ref{RD0}$). 
$$
\vCent
 {\setlength{\unitlength}{1mm}
  \begin{picture}(105,-100)
     \put(+52,-55.5){   \includegraphics[width=2.4in]{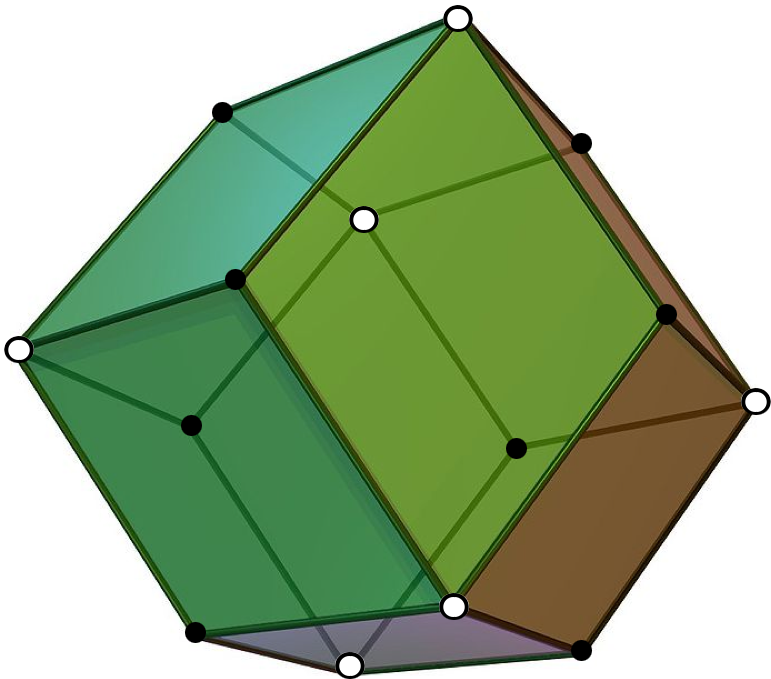}}       
  \end{picture}}
$$
\begin{equation}
\vCent
 {\setlength{\unitlength}{1mm}
  \begin{picture}(105,-100)
     \put(-21,-44.5){   \includegraphics[width=2.4in]{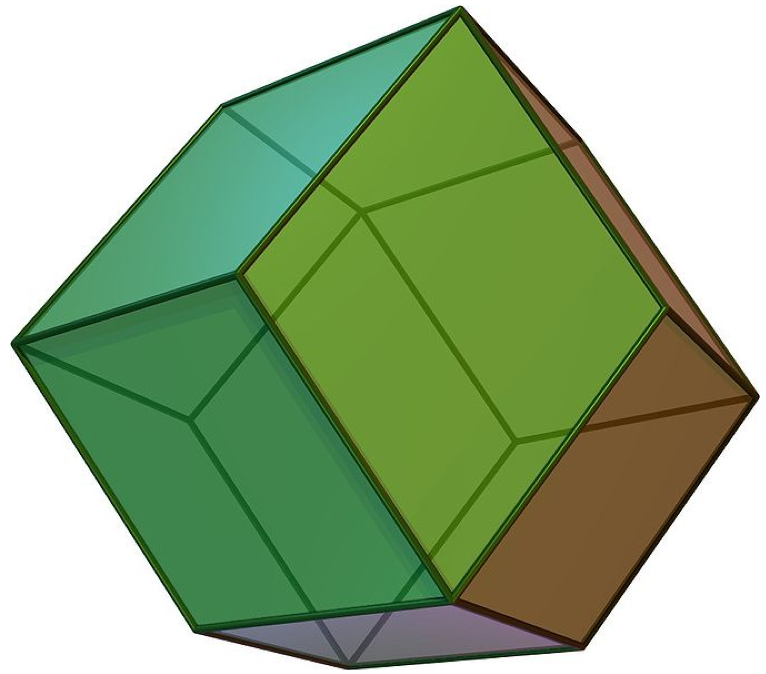}}       
  \end{picture}}
  \label{RD0}
\end{equation}
\vskip1.9in\begin{center}RD without and RD with bosonic and fermionic nodes added$~~~~~~~~~$ \end{center} 

 $~~~$Because it is a zonohedron with strictly four sided faces, we may at first expect the RD can represent an adinkra hologram that encodes a supersymmetric system of equations. Counting the white and black vertices, however, shows that we have an unequal number of fermions and bosons (see the rightmost  RD in Figure~$\ref{RD0}$). This might not seem to be an issue, but simple calculations of $\{  {\rm D}{}_{{}_{\rm I}} , {\rm D}{}_{{}_{\rm J}} \}$ proves that the graph does not follow SUSY Algebra. To show this explicitly, we must introduce the R-matrices and L-Matrices
      
 $~~~$The equations encoded into the RD valise adinkra-like graph can be represented as matrices, designated R-matrices and L-matrices, first introduced in \cite{Genomics}. The equations for the bosonic nodes with the super-differential operator, ${\rm D}{}_{{}_1}$ (the red lines), are: \\ 
\begin{center}
${\rm D}{}_{{}_1} \Phi_1 = -\dot{\imath} \Psi_3 \ \ 
{\rm D}{}_{{}_1} \Phi_2 = \dot{\imath} \Psi_5  \ \ 
{\rm D}{}_{{}_1} \Phi_3 = \dot{\imath} \Psi_1$
\be \ \ {\rm D}{}_{{}_1} \Phi_4 =\dot{\imath} \Psi_8 \ \ 
{\rm D}{}_{{}_1} \Phi_5 =\dot{\imath}\Psi_4, \ \ 
{\rm D}{}_{{}_1} \Phi_6 = -\dot{\imath} \Psi_6 \ee
\end{center}

     \vskip0.1in
\begin{equation}
\vCent
 {\setlength{\unitlength}{1mm}
  \begin{picture}(105,-100)
     \put(-0,-22.5){   \includegraphics[width=4in]{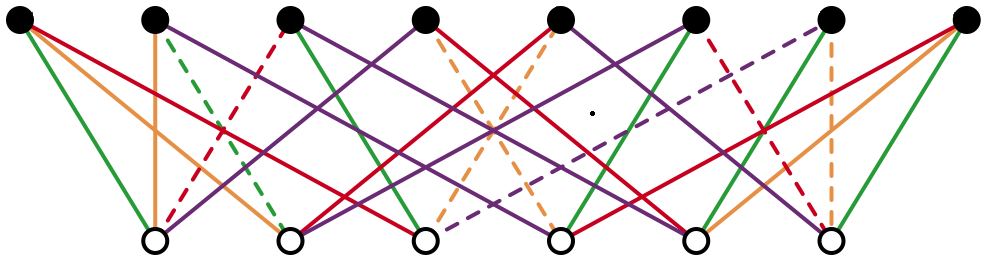}}       
  \end{picture}}
  \label{f:ValiseRD}
\nonumber
\end{equation}
\vskip0.7in\begin{center}RD Valise graph\end{center} 

This can similarly be written as:

$~~~~{\rm D}{}_1 \bold{\Phi} =\dot{\imath} ({\rm L}_1)\bold{\Psi}$

where $\bold{\Phi} =\left( \begin{array}{c}\Phi_1 \\ \Phi_2 \\ \Phi_3 \\ \Phi_4 \\ \Phi_5 \\ \Phi_6 \end{array} \right)$,  $\bold{\Psi} =\left( \begin{array}{c}\Psi_1 \\ \Psi_2 \\ \Psi_3 \\ \Psi_4 \\ \Psi_5 \\ \Psi_6\\\Psi_7 \\ \Psi_8 \end{array} \right)$ and ${\rm L}_1= \left( \begin{array}{cccccccc}
0 & 0 & -1 & 0 & 0 & 0 &0 &0\\ 
0 & 0 & 0 & 0 & 1 & 0 &0 &0\\ 
1 & 0 & 0 & 0 & 0 & 0 &0 &0\\ 
0 & 0 & 0 & 0 & 0 & 0 &0 &1\\ 
0 & 0 & 0 & 1 & 0 & 0 &0 &0\\ 
0 & 0 & 0 & 0 & 0 & -1 &0 &0\\ 
\end{array} \right)$. \\

In general, for the I-th super-differential operator, 

\be{\rm D}{}_{{}_{\rm I}} \bold{\Phi} =\dot{\imath} ({\rm L}_\rI)\bold{\Psi}\ee

for some matrix, ${\rm L}_\rI$.

The other L-Matrices for the RD (Figure ~$\ref{f:ValiseRD}$) are below:\\

\begin{center}
 ${\rm L}_2 =\left( \begin{array}{cccccccc}
0 & 0 & 0 & 1 & 0 & 0 &0 &0\\ 
0 & 0 & 0 & 0 & 0 & 1 &0 &0\\ 
0 & 0 & 0 & 0 & 0 & 0 &-1 &0\\ 
0 & 1 & 0 & 0 & 0 & 0 &0 &0\\ 
0 & 0 & 1 & 0 & 0 & 0 &0 &0\\ 
0 & 0 & 0 & 0 & 1 & 0 &0 &0\end{array} \right)$,\\
\end{center}
\be\label{eq:L}
  {\rm L}_3= \left( \begin{array}{cccccccc}
1 & 0 & 0 & 0 & 0 & 0 &0 &0\\ 
0 & -1 & 0 & 0 & 0 & 0 &0 &0\\ 
0 & 0 & 1 & 0 & 0 & 0 &0 &0\\ 
0 & 0 & 0 & 0 & 0 & 1 &0 &0\\ 
0 & 0 & 0 & 0 & 0 & 0 &1 &0\\ 
0 & 0 & 0 & 0 & 0 & 0 &0 &1\\ 
\end{array} \right), \ \
 {\rm L}_4 = \left( \begin{array}{cccccccc}
0 & 1 & 0 & 0 & 0 & 0 &0 &0\\ 
1 & 0 & 0 & 0 & 0 & 0 &0 &0\\ 
0 & 0 & 0 & 0 & -1 & 0 &0 &0\\ 
0 & 0 & 0 & -1 & 0 & 0 &0 &0\\ 
0 & 0 & 0 & 0 & 0 & 0 &0 &1\\ 
0 & 0 & 0 & 0 & 0 & 0 &-1 &0\\ 
\end{array} \right)
 \ee

The R-matrices similarly represent the relations between $\Tilde \Phi$ and $\Tilde \Psi$ for a degree flip of the valise adinkra-candidate graph\footnote{The degree flipped version of a valise adinkra graph is obtained by drawing a horizontal line \newline $~~~~~~$
 through the valise and the reflecting vertically about this line.  Since this necessarily implies
  \newline $~~~~~~$ different engineering dimensions for the fields associated with the vertices, we use
 a differ-  \newline $~~~~~~$
 ent notation for them.}:
 
$ {\rm D}{}_{{}_{\rm I}} \bold{\Tilde {\Psi}} =\dot{\imath} ({\rm R}_\rI)\bold{\Tilde {\Phi}}$

 It turns out this relationship implies the ${\rm R}_\rI$-th matrix is a transpose of the ${\rm L}_\rI$-th matrix. Therefore, the R-matrices for the valise RD are:\\

\begin{center} 
 ${\rm R}_1=  \left( \begin{array}{cccccccc}
0 & 0 & 1& 0 & 0 & 0 \\ 
0 & 0 & 0 & 0 & 0 & 0 \\ 
-1 & 0 & 0 & 0 & 0 & 0\\ 
0 & 0 & 0 & 0 & 1 & 0 \\ 
0 & 1 & 0 & 0 & 0 & 0 \\ 
0& 0 & 0 & 0 & 0 & -1 \\ 
0 & 0 & 0 & 0 & 0 & 0 \\ 
0 & 0 & 0 & 1 & 0 & 0 \\ 
\end{array} \right)$, \ \
 ${\rm R}_2 = \left( \begin{array}{cccccccc}
0 & 0 & 0 & 0 & 0 & 0 \\ 
0 & 0 & 0 & 1 & 0 & 0 \\ 
0 & 0 & 0 & 0 & 1 & 0\\ 
1 & 0 & 0 & 0 & 0 & 0 \\ 
0 & 0 & 0 & 0 & 0 & 1 \\ 
0 & 1 & 0 & 0 & 0 & 0 \\ 
0 & 0 & -1 & 0 & 0 & 0 \\ 
0 & 0 & 0 & 0 & 0 & 0 \\ 
\end{array} \right)$
\end{center}

  \be\label{eq:R}
  {\rm R}_3= \left( \begin{array}{cccccccc}
1 & 0 & 0 & 0 & 0 & 0 \\ 
0 & -1 & 0 & 0 & 0 & 0 \\ 
0 & 0 & 1 & 0 & 0 & 0\\ 
0 & 0 & 0 & 0 & 0 & 0 \\ 
0 & 0 & 0 & 0 & 0 & 0 \\ 
0 & 0 & 0 & 1 & 0 & 0 \\ 
0 & 0 & 0 & 0 & 1 & 0 \\ 
0 & 0 & 0 & 0 & 0 & 1 \\ 
\end{array} \right), \ \
{\rm R}_4 = \left( \begin{array}{cccccccc}
0 & 1 & 0 & 0 & 0 & 0 \\ 
1 & 0 & 0 & 0 & 0 & 0 \\ 
0 & 0 & 0 & 0 & 0 & 0\\ 
0 & 0 & 0 & -1 & 0 & 0 \\ 
0 & 0 & -1 & 0 & 0 & 0 \\ 
0 & 0 & 0 & 0 & 0 & 0 \\ 
0 & 0 & 0 & 0 & 0 & -1 \\ 
0 & 0 & 0 & 0 & 1 & 0 \\ 
\end{array} \right)\ee

Because the R-matrices and L-matrices are matrix versions of the super-differential operators, they realize the Garden Algebra if and only if the super-differential counterparts do as well. Therefore, in order to determine if the graphs are adinkras, we must simply check if the matrices  form a representation of the previous Garden Algebra rules:
\be\label{e:EQ2}
 (\,{\rm L}_\rI\,)_i{}^\hj\>(\,{\rm R}_\rJ\,)_\hj{}^k + (\,{\rm L}_\rJ\,)_i{}^\hj\>(\,{\rm R}_\rI\,)_\hj{}^k
   ~=~ 2 \,  \delta_{ \rI \rJ } \, \delta_i {}^k~~, 
 \ee
 \be\label{e:EQ3}
 (\,{\rm R}_\rJ\,)_\hi{}^j\>(\, {\rm L}_\rI\,)_j{}^\hk + (\,{\rm R}_\rI\,)_\hi{}^j\>(\,{\rm L}_\rJ\,)_j{}^\hk
  ~=~ 2\,\delta_{\rI\rJ}\,\delta_\hi{}^\hk~~.
\ee

The R-matrices and L-matrices in (5) and (6) do not, however, realize the Garden Algebra rules above.  Since the ratio of fermions to bosons is not equal, the R-matrices and L-matrices are not square, which then implies that the above equations create different sized matrices (as well as 1's and 0's in non-diagnal elements of the corresponding matrices). It is easy to generalize this result: \emph{any} shape with more fermion nodes than boson nodes, or vice versa, cannot be an adinkra that
describes an off-shell supermultiplet, simply because the R-matrices and L-matrices matrices would not be square.

One might try to use a double of the RD to fix this problem.  This is shown in ($\ref{f:LiftRD}$) below.  The adinkra to the left in this image may be regarded as the original RD shown in
$\ref{RD0}$.  The graph to the right is then a Klein-flipped and engineering dimensional lowered
version of the original RD-adinkra.  A fifth SUSY charge denoted by the light blue links is then
introduced in an effort to construct an adinkra for an off-shell SUSY system. This will not create an adinkra, however, because the second row of nodes will not connect to all five super-differential operators. Hence, there exists a fixed I such that:
\be \label{NOTeq}
{\rm D}{}_{\rI} \Phi_i \neq \pm \dot{\imath}\Psi_j
\ee
for some two nodes $\Phi_i$, and $\Psi_j$ for any choices of the values of $i$ and $j$, which implies that the super-differential operators do not preserve the SUSY algebra.\footnote{This ``lifting" technique, however, can be used on any zonohedron (even ones which represent \\ $~~~~~~$
adinkras). It is trivial to prove higher dimensional shapes, such as any hypercube, would pre-
 \\ $~~~~~~$
serve the Garden Algebra of the lower dimensional adinkra.  The smallest dimensional shape 
 \\ $~~~~~~$
may be called an ``Adinkra Generator" analogous to the generators in group theory. For exam-
 \\ $~~~~~~$
ple, the diamond adinkra is the generator for the hypercubic adinkras.}
   \vskip0.7in
\begin{equation}
\vCent
{\setlength{\unitlength}{1mm}
  \begin{picture}(105,-100)
     \put(15,-22.5){   \includegraphics[width=3in]{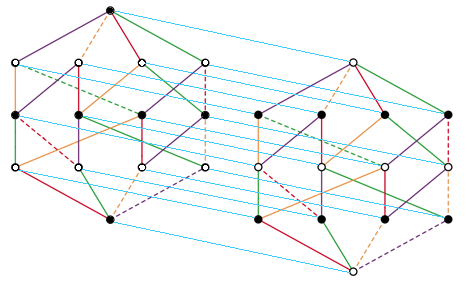}}       
  \end{picture}}
  \label{f:LiftRD}
\end{equation}
\vskip0.7in\begin{center}A double of the RD as an adinkra-like graph. Note blue line ``dashings" were removed for clarity\end{center} 

\section{The Rhombic Icosahedron: Another Zonohedron-Adinkra Candidate}

 $~~~~$ The previous case study implies we have two important restrictions in order for a zonahedron to yield an adinkra candidate that describes and off-shell SUSY multiplet,
 it must have:\newline
$~~~~~$-strictly four sided faces, and

-equal numbers of fermions and bosons.

Although it might be considered a sad ending to our search if there was no minimal dimensional zonohedron which could create an adinkra, in fact, this is in accord with the ``folk wisdom'' that
has been seen all previous studies.

Clearly, two colors admit only the trivial the bow-tie and diamond adinkras. With three colors, the search for eligible shapes broadens significantly, so one could hope that with this number of
colors, at least, one might find  another adinkra shape, besides a cube.  Within the class of zonohedra,
this search leads to the the rhombic icosahedron (to be designated as  ``RI". See Figure~$\ref{f:IcosRD}$). 
   \vskip0.7in
$$
\vCent
{\setlength{\unitlength}{1mm}
  \begin{picture}(105,-100)
     \put(-32,-59.5){   \includegraphics[width=2.8in]{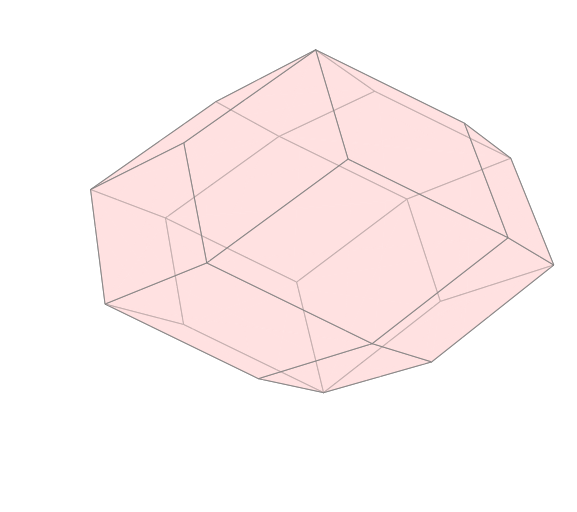}}       
  \end{picture}}
  \label{f:IcosRD}
$$
\begin{equation}
\vCent
{\setlength{\unitlength}{1mm}
  \begin{picture}(105,-100)
     \put(44,-51){   \includegraphics[width=2.8in]{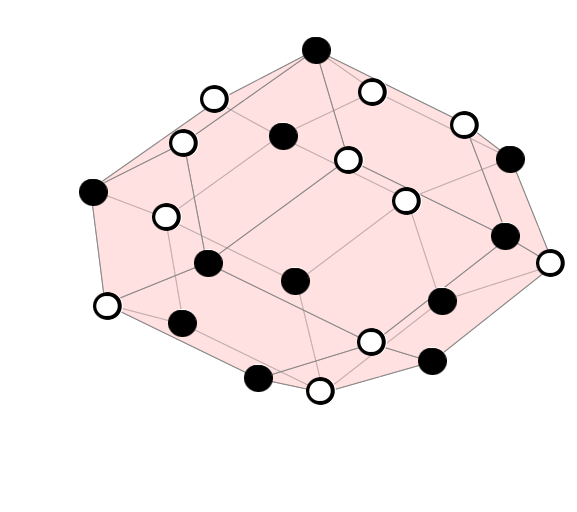}}       
  \end{picture}}
  \label{f:IcosRD}
\end{equation}
\vskip1.7in\begin{center}The RI with \& without boson and fermion nodes\end{center} 

A simple study will reveal that the RI not only has strictly four sided faces, but also equal numbers of fermions and bosons.  To determine whether it satisfied the Garden Algebra, we checked whether the R-matrices and L-matrices of the RI satisfied ($\ref{e:EQ2}$) and ($\ref{e:EQ3}$).\\\\
The L-matrices were found to be:\newline\newline
\be {\rm L}_1=\left(
\begin{array}{ccccccccccc}
 -1 & 0 & 0 & 0 & 0 & 0 & 0 & 0 & 0 & 0 & 0 \\
 0 & 1 & 0 & 0 & 0 & 0 & 0 & 0 & 0 & 0 & 0 \\
 0 & 0 & 0 & 0 & 0 & 0 & 0 & 0 & 0 & 0 & 0 \\
 0 & 0 & 0 & 0 & 0 & 0 & 0 & 0 & 0 & 0 & 0 \\
 0 & 0 & 0 & -1 & 0 & 0 & 0 & 0 & 0 & 0 & 0 \\
 0 & 0 & -1 & 0 & 0 & 0 & 0 & 0 & 0 & 0 & 0 \\
 0 & 0 & 0 & 0 & 0 & 0 & 1 & 0 & 0 & 0 & 0 \\
 0 & 0 & 0 & 0 & 0 & 0 & 0 & 0 & 0 & 0 & 0 \\
 0 & 0 & 0 & 0 & 0 & 0 & 0 & 0 & 0 & 1 & 0 \\
 0 & 0 & 0 & 0 & 0 & 1 & 0 & 0 & 0 & 0 & 0 \\
 0 & 0 & 0 & 0 & 0 & 0 & 0 & -1 & 0 & 0 & 0
\end{array}
\right)\ee
\be {\rm L}_2=\left(
\begin{array}{ccccccccccc}
 0 & 1 & 0 & 0 & 0 & 0 & 0 & 0 & 0 & 0 & 0 \\
 1 & 0 & 0 & 0 & 0 & 0 & 0 & 0 & 0 & 0 & 0 \\
 0 & 0 & 1 & 0 & 0 & 0 & 0 & 0 & 0 & 0 & 0 \\
 0 & 0 & 0 & 0 & 0 & 0 & 0 & 0 & 0 & 0 & 0 \\
 0 & 0 & 0 & 0 & 0 & 0 & 0 & 0 & 0 & 0 & 0 \\
 0 & 0 & 0 & 0 & 0 & 0 & 0 & 0 & 0 & 0 & 0 \\
 0 & 0 & 0 & 0 & -1 & 0 & 0 & 0 & 0 & 0 & 0 \\
 0 & 0 & 0 & -1 & 0 & 0 & 0 & 0 & 0 & 0 & 0 \\
 0 & 0 & 0 & 0 & 0 & 0 & 0 & 1 & 0 & 0 & 0 \\
 0 & 0 & 0 & 0 & 0 & 0 & 0 & 0 & 0 & 0 & 1 \\
 0 & 0 & 0 & 0 & 0 & 0 & 0 & 0 & 0 & 1 & 0
\end{array}
\right)\ee
\be {\rm L}_3=\left(
\begin{array}{ccccccccccc}
 0 & 0 & 0 & 0 & 0 & 0 & 0 & 0 & 0 & 0 & 0 \\
 0 & 0 & 1 & 0 & 0 & 0 & 0 & 0 & 0 & 0 & 0 \\
 -1 & 0 & 0 & 0 & 0 & 0 & 0 & 0 & 0 & 0 & 0 \\
 0 & 0 & 0 & 0 & 1 & 0 & 0 & 0 & 0 & 0 & 0 \\
 0 & 0 & 0 & 0 & 0 & 0 & 0 & 0 & 0 & 0 & 0 \\
 0 & 1 & 0 & 0 & 0 & 0 & 0 & 0 & 0 & 0 & 0 \\
 0 & 0 & 0 & 0 & 0 & 0 & 0 & 0 & 0 & 0 & 0 \\
 0 & 0 & 0 & 0 & 0 & 0 & 0 & 0 & 1 & 0 & 0 \\
 0 & 0 & 0 & 0 & 0 & -1 & 0 & 0 & 0 & 0 & 0 \\
 0 & 0 & 0 & 0 & 0 & 0 & 0 & 0 & 0 & 1 & 0 \\
 0 & 0 & 0 & 0 & 0 & 0 & 0 & 0 & 0 & 0 & -1
\end{array}
\right)\ee
\be {\rm L}_4=\left(
\begin{array}{ccccccccccc}
 0 & 0 & 0 & 0 & 0 & 0 & 0 & 0 & 0 & 0 & 0 \\
 0 & 0 & 0 & 0 & 0 & 0 & 0 & 0 & 0 & 0 & 0 \\
 0 & 0 & 0 & 0 & 1 & 0 & 0 & 0 & 0 & 0 & 0 \\
 1 & 0 & 0 & 0 & 0 & 0 & 0 & 0 & 0 & 0 & 0 \\
 0 & 0 & 0 & 0 & 0 & 1 & 0 & 0 & 0 & 0 & 0 \\
 0 & 0 & 0 & 0 & 0 & 0 & 1 & 0 & 0 & 0 & 0 \\
 0 & 0 & 1 & 0 & 0 & 0 & 0 & 0 & 0 & 0 & 0 \\
 0 & 0 & 0 & 0 & 0 & 0 & 0 & 0 & 0 & 0 & 1 \\
 0 & 0 & 0 & 0 & 0 & 0 & 0 & 0 & 0 & 0 & 0 \\
 0 & 0 & 0 & 1 & 0 & 0 & 0 & 0 & 0 & 0 & 0 \\
 0 & 0 & 0 & 0 & 0 & 0 & 0 & 0 & 1 & 0 & 0
\end{array}
\right)\ee
\be {\rm L}_5=\left(
\begin{array}{ccccccccccc}
 0 & 0 & 0 & 1 & 0 & 0 & 0 & 0 & 0 & 0 & 0 \\
 0 & 0 & 0 & 0 & 0 & 0 & 0 & 0 & 0 & 0 & 0 \\
 0 & 0 & 0 & 0 & 0 & 0 & 0 & 0 & 0 & 0 & 0 \\
 0 & 0 & 0 & 0 & 0 & 1 & 0 & 0 & 0 & 0 & 0 \\
 -1 & 0 & 0 & 0 & 0 & 0 & 0 & 0 & 0 & 0 & 0 \\
 0 & 0 & 0 & 0 & 0 & 0 & 0 & 0 & -1 & 0 & 0 \\
 0 & 0 & 0 & 0 & 0 & 0 & 0 & 1 & 0 & 0 & 0 \\
 0 & 1 & 0 & 0 & 0 & 0 & 0 & 0 & 0 & 0 & 0 \\
 0 & 0 & 0 & 0 & 1 & 0 & 0 & 0 & 0 & 0 & 0 \\
 0 & 0 & 0 & 0 & 0 & 0 & 0 & 0 & 0 & 0 & 0 \\
 0 & 0 & 0 & 0 & 0 & 0 & 1 & 0 & 0 & 0 & 0
\end{array}
\right)\ee

The R-matrices are therefore:\newline

\be {\rm R}_1=   \left(
\begin{array}{ccccccccccc}
 -1 & 0 & 0 & 0 & 0 & 0 & 0 & 0 & 0 & 0 & 0 \\
 0 & 1 & 0 & 0 & 0 & 0 & 0 & 0 & 0 & 0 & 0 \\
 0 & 0 & 0 & 0 & 0 & -1 & 0 & 0 & 0 & 0 & 0 \\
 0 & 0 & 0 & 0 & -1 & 0 & 0 & 0 & 0 & 0 & 0 \\
 0 & 0 & 0 & 0 & 0 & 0 & 0 & 0 & 0 & 0 & 0 \\
 0 & 0 & 0 & 0 & 0 & 0 & 0 & 0 & 0 & 1 & 0 \\
 0 & 0 & 0 & 0 & 0 & 0 & 1 & 0 & 0 & 0 & 0 \\
 0 & 0 & 0 & 0 & 0 & 0 & 0 & 0 & 0 & 0 & -1 \\
 0 & 0 & 0 & 0 & 0 & 0 & 0 & 0 & 0 & 0 & 0 \\
 0 & 0 & 0 & 0 & 0 & 0 & 0 & 0 & 1 & 0 & 0 \\
 0 & 0 & 0 & 0 & 0 & 0 & 0 & 0 & 0 & 0 & 0
\end{array}
\right)\ee 
\be {\rm R}_2=\left(
\begin{array}{ccccccccccc}
 0 & 1 & 0 & 0 & 0 & 0 & 0 & 0 & 0 & 0 & 0 \\
 1 & 0 & 0 & 0 & 0 & 0 & 0 & 0 & 0 & 0 & 0 \\
 0 & 0 & 1 & 0 & 0 & 0 & 0 & 0 & 0 & 0 & 0 \\
 0 & 0 & 0 & 0 & 0 & 0 & 0 & -1 & 0 & 0 & 0 \\
 0 & 0 & 0 & 0 & 0 & 0 & -1 & 0 & 0 & 0 & 0 \\
 0 & 0 & 0 & 0 & 0 & 0 & 0 & 0 & 0 & 0 & 0 \\
 0 & 0 & 0 & 0 & 0 & 0 & 0 & 0 & 0 & 0 & 0 \\
 0 & 0 & 0 & 0 & 0 & 0 & 0 & 0 & 1 & 0 & 0 \\
 0 & 0 & 0 & 0 & 0 & 0 & 0 & 0 & 0 & 0 & 0 \\
 0 & 0 & 0 & 0 & 0 & 0 & 0 & 0 & 0 & 0 & 1 \\
 0 & 0 & 0 & 0 & 0 & 0 & 0 & 0 & 0 & 1 & 0
\end{array}
\right)\ee
\be {\rm R}_3=\left(
\begin{array}{ccccccccccc}
 0 & 0 & -1 & 0 & 0 & 0 & 0 & 0 & 0 & 0 & 0 \\
 0 & 0 & 0 & 0 & 0 & 1 & 0 & 0 & 0 & 0 & 0 \\
 0 & 1 & 0 & 0 & 0 & 0 & 0 & 0 & 0 & 0 & 0 \\
 0 & 0 & 0 & 0 & 0 & 0 & 0 & 0 & 0 & 0 & 0 \\
 0 & 0 & 0 & 1 & 0 & 0 & 0 & 0 & 0 & 0 & 0 \\
 0 & 0 & 0 & 0 & 0 & 0 & 0 & 0 & -1 & 0 & 0 \\
 0 & 0 & 0 & 0 & 0 & 0 & 0 & 0 & 0 & 0 & 0 \\
 0 & 0 & 0 & 0 & 0 & 0 & 0 & 0 & 0 & 0 & 0 \\
 0 & 0 & 0 & 0 & 0 & 0 & 0 & 1 & 0 & 0 & 0 \\
 0 & 0 & 0 & 0 & 0 & 0 & 0 & 0 & 0 & 1 & 0 \\
 0 & 0 & 0 & 0 & 0 & 0 & 0 & 0 & 0 & 0 & -1
\end{array}
\right)\ee 
\be {\rm R}_4=\left(
\begin{array}{ccccccccccc}
 0 & 0 & 0 & 1 & 0 & 0 & 0 & 0 & 0 & 0 & 0 \\
 0 & 0 & 0 & 0 & 0 & 0 & 0 & 0 & 0 & 0 & 0 \\
 0 & 0 & 0 & 0 & 0 & 0 & 1 & 0 & 0 & 0 & 0 \\
 0 & 0 & 0 & 0 & 0 & 0 & 0 & 0 & 0 & 1 & 0 \\
 0 & 0 & 1 & 0 & 0 & 0 & 0 & 0 & 0 & 0 & 0 \\
 0 & 0 & 0 & 0 & 1 & 0 & 0 & 0 & 0 & 0 & 0 \\
 0 & 0 & 0 & 0 & 0 & 1 & 0 & 0 & 0 & 0 & 0 \\
 0 & 0 & 0 & 0 & 0 & 0 & 0 & 0 & 0 & 0 & 0 \\
 0 & 0 & 0 & 0 & 0 & 0 & 0 & 0 & 0 & 0 & 1 \\
 0 & 0 & 0 & 0 & 0 & 0 & 0 & 0 & 0 & 0 & 0 \\
 0 & 0 & 0 & 0 & 0 & 0 & 0 & 1 & 0 & 0 & 0
\end{array}\right)\ee
\be {\rm R}_5=\left(
\begin{array}{ccccccccccc}
 0 & 0 & 0 & 0 & -1 & 0 & 0 & 0 & 0 & 0 & 0 \\
 0 & 0 & 0 & 0 & 0 & 0 & 0 & 1 & 0 & 0 & 0 \\
 0 & 0 & 0 & 0 & 0 & 0 & 0 & 0 & 0 & 0 & 0 \\
 1 & 0 & 0 & 0 & 0 & 0 & 0 & 0 & 0 & 0 & 0 \\
 0 & 0 & 0 & 0 & 0 & 0 & 0 & 0 & 1 & 0 & 0 \\
 0 & 0 & 0 & 1 & 0 & 0 & 0 & 0 & 0 & 0 & 0 \\
 0 & 0 & 0 & 0 & 0 & 0 & 0 & 0 & 0 & 0 & 1 \\
 0 & 0 & 0 & 0 & 0 & 0 & 1 & 0 & 0 & 0 & 0 \\
 0 & 0 & 0 & 0 & 0 & -1 & 0 & 0 & 0 & 0 & 0 \\
 0 & 0 & 0 & 0 & 0 & 0 & 0 & 0 & 0 & 0 & 0 \\
 0 & 0 & 0 & 0 & 0 & 0 & 0 & 0 & 0 & 0 & 0
\end{array}
\right)\ee\\\\
It is easy to show that these matrices \emph{do not} satisfy Garden Algebra (see Appendix B). Alike to the four-dimensional version of the RD, only two nodes connect to all the super-differential operators.

The Icosahedron is important to demonstrate the final restriction necessary for a three dimensional shape to represent an adinkra. Because we need rhombic faces for a polyhedron to represent a SUSY adinkra, the second row of the polyhedron's graph will connect to no more than three other nodes. This forces us to only consider shapes with three super-differential operators, all of which are topologically the same as a cube. 

This not only rules out all other possible three dimensional shapes, but it allows us to rule out all non-hypercubic zonotopes as adinkra candidates. For example, just as a three color zonohedron must have faces which are topologically the same as diamond adinkras, four color zonotope must therefore have ``faces" which are topologically the same as cubes. Also, just as before, every adinkra candidate must have the same number of fermions and bosons. Lastly, the second row of a four dimensional zonotope cannot connect to more than four other nodes, which then implies the tesseract is the only valid adinkra.

\section{Conclusion}

$~~~~$ We can now summarize what was found. Firstly, there is valid reason to suspect that non-hypercube adinkras exist. To find, at least initially, what shapes could be adinkras, we search for strictly four sided faced zonohedra. Even under this restriction, we must also search for zonohedra with equal numbers of fermions and bosons on its vertices, before we actually check and see if the adinkra satisfies Garden Algebra. Lastly, every super-differential operator must connect to each node. Because the second row cannot have more than three super-differential operators connected to it, when represented by a polyhedron, we must have no more than three differential operators. This forces us to conclude that the cube is the only valid polyhedron that can to create adinkras. Higher dimensional analogues can then be created, which restrict the zonotopes to only be hypercubes.

There is one point in our discussion that could bare further scrutiny.  The choices of ``dashedness''
i.\ e.\ the signs of the non-vanishing entries in the L-matrices and R-matrices are somewhat arbitrary.  In the case of the RD, we use a ``trick'' to make the assignment.  The RD may be regarded as the tesseract where two bosonic nodes separated by the greatest Hamming distance have been deleted.
This is how we fixed the signs for the RD.  For the RI there is a similar minimal ansatz that can
be used.  We believe that none of our major conclusions are affected by other choices of 
``dashedness.''

${~~~}$ \newline
${~~~~~}$``{\it {Colors answer feeling in man; shapes answer thought; \newline
${~~~~~~}$and motion answers will.}}" \newline
\newline $~~~~~~~$ -- John Sterling
\newline ${~~~}$ \newline

\noindent
{\Large\bf Acknowledgments} 
This research was supported in part by the endowment of the John S.~Toll Professorship, the University of Maryland Center for String \& Particle Theory, National Science Foundation Grant PHY-0354401.  SJG  and KS offer additional gratitude to the M.\ L.\ K. Visiting Professorship and to the M.\ I.\ T.\ Center for Theoretical Physics for support and hospitality facilitating part of this work.  This work is also supported 
by U.S. Department of Energy (D.O.E.) under cooperative agreement DEFG0205ER41360.

\newpage
\section{Appendix A}

$~~~~~$ Here, we show that the the R-matrices and L-matrices of the rhombic dodecahedron don't satisfy the Garden Algebra, and therefore cannot lead to adinkras that describe off-shell SUSY 
systems.  The computations of the L-matrices and R-matrices are below.
\be \label{L_1 R_1} {\rm L}_1*{\rm R}_1 = \left(
\begin{array}{cccccc}
 1 & 0 & 0 & 0 & 0 & 0 \\
 0 & 1 & 0 & 0 & 0 & 0 \\
 0 & 0 & 1 & 0 & 0 & 0 \\
 0 & 0 & 0 & 1 & 0 & 0 \\
 0 & 0 & 0 & 0 & 1 & 0 \\
 0 & 0 & 0 & 0 & 0 & 1
\end{array}
\right)\ee 
\be \label{L_1R_2L_2R_1} {\rm L}_1*{\rm R}_2 + {\rm L}_2*{\rm R}_1 =\left(
\begin{array}{cccccc}
 0 & 0 & 0 & 0 & 0 & 0 \\
 0 & 0 & 0 & 0 & 0 & 0 \\
 0 & 0 & 0 & 0 & 0 & 0 \\
 0 & 0 & 0 & 0 & 0 & 0 \\
 0 & 0 & 0 & 0 & 0 & 0 \\
 0 & 0 & 0 & 0 & 0 & 0
\end{array}
\right)\ee 
\be\label{L_1R_3L_3R_1} {\rm L}_1*{\rm R}_3+{\rm L}_3*{\rm R}_1=\left(
\begin{array}{cccccc}
 0 & 0 & 0 & 0 & 0 & 0 \\
 0 & 0 & 0 & 0 & 0 & 0 \\
 0 & 0 & 0 & 0 & 0 & 0 \\
 0 & 0 & 0 & 0 & 0 & 0 \\
 0 & 0 & 0 & 0 & 0 & 0 \\
 0 & 0 & 0 & 0 & 0 & 0
\end{array}
\right)\ee 
\be \label{L_1R_3L_3R_1}
{\rm L}_1 *{\rm R}_4 + {\rm L}_4 *{\rm R}_1=\left(
\begin{array}{cccccc}
 0 & 0 & 0 & 0 & 0 & 0 \\
 0 & 0 & 0 & 0 & 0 & 0 \\
 0 & 0 & 0 & 0 & 0 & 0 \\
 0 & 0 & 0 & 0 & 0 & 0 \\
 0 & 0 & 0 & 0 & 0 & 0 \\
 0 & 0 & 0 & 0 & 0 & 0
\end{array}
\right)\ee 
\be\label{L_2R_2} {\rm L}_2*{\rm R}_2=\left(
\begin{array}{cccccc}
 1 & 0 & 0 & 0 & 0 & 0 \\
 0 & 1 & 0 & 0 & 0 & 0 \\
 0 & 0 & 1 & 0 & 0 & 0 \\
 0 & 0 & 0 & 1 & 0 & 0 \\
 0 & 0 & 0 & 0 & 1 & 0 \\
 0 & 0 & 0 & 0 & 0 & 1
\end{array}
\right)\ee 
\be \label{L_2R_3L_3R_2}{\rm L}_2 *{\rm R}_3 + {\rm L}_3 *{\rm R}_2=\left(
\begin{array}{cccccc}
 0 & 0 & 0 & 0 & 0 & 0 \\
 0 & 0 & 0 & 0 & 0 & 0 \\
 0 & 0 & 0 & 0 & 0 & 0 \\
 0 & 0 & 0 & 0 & 0 & 0 \\
 0 & 0 & 0 & 0 & 0 & 0 \\
 0 & 0 & 0 & 0 & 0 & 0
\end{array}
\right)\ee 
\be  \label{L_2R_4L_4R_2} {\rm L}_2*{\rm R}_4+{\rm L}_4*{\rm R}_2=\left(
\begin{array}{cccccc}
 0 & 0 & 0 & 0 & 0 & 0 \\
 0 & 0 & 0 & 0 & 0 & 0 \\
 0 & 0 & 0 & 0 & 0 & 0 \\
 0 & 0 & 0 & 0 & 0 & 0 \\
 0 & 0 & 0 & 0 & 0 & 0 \\
 0 & 0 & 0 & 0 & 0 & 0
\end{array}
\right)\ee 
\be \label{L_3R_3} {\rm L}_3 *{\rm R}_3=\left(
\begin{array}{cccccc}
 1 & 0 & 0 & 0 & 0 & 0 \\
 0 & 1 & 0 & 0 & 0 & 0 \\
 0 & 0 & 1 & 0 & 0 & 0 \\
 0 & 0 & 0 & 1 & 0 & 0 \\
 0 & 0 & 0 & 0 & 1 & 0 \\
 0 & 0 & 0 & 0 & 0 & 1
\end{array}
\right)\ee 
\be\label{L_3R_4L_4R_3} {\rm L}_3*{\rm R}_4+{\rm L}_4*{\rm R}_3= \left(
\begin{array}{cccccc}
 0 & 0 & 0 & 0 & 0 & 0 \\
 0 & 0 & 0 & 0 & 0 & 0 \\
 0 & 0 & 0 & 0 & 0 & 0 \\
 0 & 0 & 0 & 0 & 0 & 0 \\
 0 & 0 & 0 & 0 & 0 & 0 \\
 0 & 0 & 0 & 0 & 0 & 0
\end{array}
\right)\ee 
\be\label{L_4R_4} {\rm L}_4*{\rm R}_4 = \left(
\begin{array}{cccccc}
 1 & 0 & 0 & 0 & 0 & 0 \\
 0 & 1 & 0 & 0 & 0 & 0 \\
 0 & 0 & 1 & 0 & 0 & 0 \\
 0 & 0 & 0 & 1 & 0 & 0 \\
 0 & 0 & 0 & 0 & 1 & 0 \\
 0 & 0 & 0 & 0 & 0 & 1
\end{array}
\right)\ee 

We can see the matrices satisfy the first equation in (\ref{eq:EQ1}). 

The rest of the matrices are computed below, as a counter-example to the claim that these could satisfy the equation in (\ref{e:EQ3}), thereby proving the rhombic dodecahedron cannot satisfy Garden Algebra:

\be \label{L_1R_1} {\rm R}_1*{\rm L}_1=\left(
\begin{array}{cccccccc}
 1 & 0 & 0 & 0 & 0 & 0 & 0 & 0 \\
 0 & 0 & 0 & 0 & 0 & 0 & 0 & 0 \\
 0 & 0 & 1 & 0 & 0 & 0 & 0 & 0 \\
 0 & 0 & 0 & 1 & 0 & 0 & 0 & 0 \\
 0 & 0 & 0 & 0 & 1 & 0 & 0 & 0 \\
 0 & 0 & 0 & 0 & 0 & 1 & 0 & 0 \\
 0 & 0 & 0 & 0 & 0 & 0 & 0 & 0 \\
 0 & 0 & 0 & 0 & 0 & 0 & 0 & 1
\end{array}
\right)\ee 
\be\label{R_1L_2R_2L_1} {\rm R}_1*{\rm L}_2 + {\rm R}_2*{\rm L}_1=\left(
\begin{array}{cccccccc}
 0 & 0 & 0 & 0 & 0 & 0 & -1 & 0 \\
 0 & 0 & 0 & 0 & 0 & 0 & 0 & 1 \\
 0 & 0 & 0 & 0 & 0 & 0 & 0 & 0 \\
 0 & 0 & 0 & 0 & 0 & 0 & 0 & 0 \\
 0 & 0 & 0 & 0 & 0 & 0 & 0 & 0 \\
 0 & 0 & 0 & 0 & 0 & 0 & 0 & 0 \\
 -1 & 0 & 0 & 0 & 0 & 0 & 0 & 0 \\
 0 & 1 & 0 & 0 & 0 & 0 & 0 & 0
\end{array}
\right)\ee 
\be\label{R_1L_3R_3L_1} {\rm R}_1*{\rm L}_3+{\rm R}_3 *{\rm L}_1= \left(
\begin{array}{cccccccc}
 0 & 0 & 0 & 0 & 0 & 0 & 0 & 0 \\
 0 & 0 & 0 & 0 & -1 & 0 & 0 & 0 \\
 0 & 0 & 0 & 0 & 0 & 0 & 0 & 0 \\
 0 & 0 & 0 & 0 & 0 & 0 & 1 & 0 \\
 0 & -1 & 0 & 0 & 0 & 0 & 0 & 0 \\
 0 & 0 & 0 & 0 & 0 & 0 & 0 & 0 \\
 0 & 0 & 0 & 1 & 0 & 0 & 0 & 0 \\
 0 & 0 & 0 & 0 & 0 & 0 & 0 & 0
\end{array}
\right)\ee 
\be  \label{R_1L_4R_4L_1}{\rm R}_1 *{\rm L}_4 + {\rm R}_4 *{\rm L}_1=\left(
\begin{array}{cccccccc}
 0 & 0 & 0 & 0 & 0 & 0 & 0 & 0 \\
 0 & 0 & -1 & 0 & 0 & 0 & 0 & 0 \\
 0 & -1 & 0 & 0 & 0 & 0 & 0 & 0 \\
 0 & 0 & 0 & 0 & 0 & 0 & 0 & 0 \\
 0 & 0 & 0 & 0 & 0 & 0 & 0 & 0 \\
 0 & 0 & 0 & 0 & 0 & 0 & 1 & 0 \\
 0 & 0 & 0 & 0 & 0 & 1 & 0 & 0 \\
 0 & 0 & 0 & 0 & 0 & 0 & 0 & 0
\end{array}
\right)\ee 
\be\label{R_2L_2} {\rm R}_2*{\rm L}_2=\left(
\begin{array}{cccccccc}
 0 & 0 & 0 & 0 & 0 & 0 & 0 & 0 \\
 0 & 1 & 0 & 0 & 0 & 0 & 0 & 0 \\
 0 & 0 & 1 & 0 & 0 & 0 & 0 & 0 \\
 0 & 0 & 0 & 1 & 0 & 0 & 0 & 0 \\
 0 & 0 & 0 & 0 & 1 & 0 & 0 & 0 \\
 0 & 0 & 0 & 0 & 0 & 1 & 0 & 0 \\
 0 & 0 & 0 & 0 & 0 & 0 & 1 & 0 \\
 0 & 0 & 0 & 0 & 0 & 0 & 0 & 0
\end{array}
\right)\ee 
\be\label{R_2L_3R_3L_2} {\rm R}_2*{\rm L}_3 + {\rm R}_3*{\rm L}_2=\left(
\begin{array}{cccccccc}
 0 & 0 & 0 & 1 & 0 & 0 & 0 & 0 \\
 0 & 0 & 0 & 0 & 0 & 0 & 0 & 0 \\
 0 & 0 & 0 & 0 & 0 & 0 & 0 & 0 \\
 1 & 0 & 0 & 0 & 0 & 0 & 0 & 0 \\
 0 & 0 & 0 & 0 & 0 & 0 & 0 & 1 \\
 0 & 0 & 0 & 0 & 0 & 0 & 0 & 0 \\
 0 & 0 & 0 & 0 & 0 & 0 & 0 & 0 \\
 0 & 0 & 0 & 0 & 1 & 0 & 0 & 0
\end{array}
\right)\ee 
\be \label{R_2L_4R_4L_2} {\rm R}_2*{\rm L}_4 + {\rm R}_4*{\rm L}_2=\left(
\begin{array}{cccccccc}
 0 & 0 & 0 & 0 & 0 & 1 & 0 & 0 \\
 0 & 0 & 0 & 0 & 0 & 0 & 0 & 0 \\
 0 & 0 & 0 & 0 & 0 & 0 & 0 & 1 \\
 0 & 0 & 0 & 0 & 0 & 0 & 0 & 0 \\
 0 & 0 & 0 & 0 & 0 & 0 & 0 & 0 \\
 1 & 0 & 0 & 0 & 0 & 0 & 0 & 0 \\
 0 & 0 & 0 & 0 & 0 & 0 & 0 & 0 \\
 0 & 0 & 1 & 0 & 0 & 0 & 0 & 0
\end{array}
\right)\ee 

\be \label{R_3L_3} {\rm R}_3*{\rm L}_3=\left(
\begin{array}{cccccccc}
 1 & 0 & 0 & 0 & 0 & 0 & 0 & 0 \\
 0 & 1 & 0 & 0 & 0 & 0 & 0 & 0 \\
 0 & 0 & 1 & 0 & 0 & 0 & 0 & 0 \\
 0 & 0 & 0 & 0 & 0 & 0 & 0 & 0 \\
 0 & 0 & 0 & 0 & 0 & 0 & 0 & 0 \\
 0 & 0 & 0 & 0 & 0 & 1 & 0 & 0 \\
 0 & 0 & 0 & 0 & 0 & 0 & 1 & 0 \\
 0 & 0 & 0 & 0 & 0 & 0 & 0 & 1
\end{array}
\right)\ee 
\be  \label{R_3L_4R_4L_3}{\rm R}_3*{\rm L}_4+{\rm R}_4*{\rm L}_3=\left(
\begin{array}{cccccccc}
 0 & 0 & 0 & 0 & 0 & 0 & 0 & 0 \\
 0 & 0 & 0 & 0 & 0 & 0 & 0 & 0 \\
 0 & 0 & 0 & 0 & -1 & 0 & 0 & 0 \\
 0 & 0 & 0 & 0 & 0 & -1 & 0 & 0 \\
 0 & 0 & -1 & 0 & 0 & 0 & 0 & 0 \\
 0 & 0 & 0 & -1 & 0 & 0 & 0 & 0 \\
 0 & 0 & 0 & 0 & 0 & 0 & 0 & 0 \\
 0 & 0 & 0 & 0 & 0 & 0 & 0 & 0
\end{array}
\right)\ee 
\be  \label{R4L4}{\rm R}_4*{\rm L}_4=\left(
\begin{array}{cccccccc}
 1 & 0 & 0 & 0 & 0 & 0 & 0 & 0 \\
 0 & 1 & 0 & 0 & 0 & 0 & 0 & 0 \\
 0 & 0 & 0 & 0 & 0 & 0 & 0 & 0 \\
 0 & 0 & 0 & 1 & 0 & 0 & 0 & 0 \\
 0 & 0 & 0 & 0 & 1 & 0 & 0 & 0 \\
 0 & 0 & 0 & 0 & 0 & 0 & 0 & 0 \\
 0 & 0 & 0 & 0 & 0 & 0 & 1 & 0 \\
 0 & 0 & 0 & 0 & 0 & 0 & 0 & 1
\end{array}
\right)\ee 

\section{Appendix B}

Similar to the discussion in Appendix A, we compute the relevant R-matrices and L-matrices for the Rhombic Icosahedron. These matrices clearly do not satisfy (\ref{e:EQ2}), or (\ref{e:EQ3}).
The explicit calculations follow below: \\    \\
    
\be \label{L_1R_1} {\rm L}_1*{\rm R}_1 = \left(
\begin{array}{ccccccccccc}
 1 & 0 & 0 & 0 & 0 & 0 & 0 & 0 & 0 & 0 & 0 \\
 0 & 1 & 0 & 0 & 0 & 0 & 0 & 0 & 0 & 0 & 0 \\
 0 & 0 & 0 & 0 & 0 & 0 & 0 & 0 & 0 & 0 & 0 \\
 0 & 0 & 0 & 0 & 0 & 0 & 0 & 0 & 0 & 0 & 0 \\
 0 & 0 & 0 & 0 & 1 & 0 & 0 & 0 & 0 & 0 & 0 \\
 0 & 0 & 0 & 0 & 0 & 1 & 0 & 0 & 0 & 0 & 0 \\
 0 & 0 & 0 & 0 & 0 & 0 & 1 & 0 & 0 & 0 & 0 \\
 0 & 0 & 0 & 0 & 0 & 0 & 0 & 0 & 0 & 0 & 0 \\
 0 & 0 & 0 & 0 & 0 & 0 & 0 & 0 & 1 & 0 & 0 \\
 0 & 0 & 0 & 0 & 0 & 0 & 0 & 0 & 0 & 1 & 0 \\
 0 & 0 & 0 & 0 & 0 & 0 & 0 & 0 & 0 & 0 & 1
\end{array}
\right)\ee 
\be \label{L_1R_2L_2R_1} {\rm L}_1*{\rm R}_2 + {\rm L}_2*{\rm R}_1 =\left(
\begin{array}{ccccccccccc}
 0 & 0 & 0 & 0 & 0 & 0 & 0 & 0 & 0 & 0 & 0 \\
 0 & 0 & 0 & 0 & 0 & 0 & 0 & 0 & 0 & 0 & 0 \\
 0 & 0 & 0 & 0 & 0 & -1 & 0 & 0 & 0 & 0 & 0 \\
 0 & 0 & 0 & 0 & 0 & 0 & 0 & 0 & 0 & 0 & 0 \\
 0 & 0 & 0 & 0 & 0 & 0 & 0 & 1 & 0 & 0 & 0 \\
 0 & 0 & -1 & 0 & 0 & 0 & 0 & 0 & 0 & 0 & 0 \\
 0 & 0 & 0 & 0 & 0 & 0 & 0 & 0 & 0 & 0 & 0 \\
 0 & 0 & 0 & 0 & 1 & 0 & 0 & 0 & 0 & 0 & 0 \\
 0 & 0 & 0 & 0 & 0 & 0 & 0 & 0 & 0 & 0 & 0 \\
 0 & 0 & 0 & 0 & 0 & 0 & 0 & 0 & 0 & 0 & 0 \\
 0 & 0 & 0 & 0 & 0 & 0 & 0 & 0 & 0 & 0 & 0
\end{array}
\right)\ee 
\be\label{L_1R_3L_3R_1} {\rm L}_1*{\rm R}_3+{\rm L}_3*{\rm R}_1=\left(
\begin{array}{ccccccccccc}
 0 & 0 & 1 & 0 & 0 & 0 & 0 & 0 & 0 & 0 & 0 \\
 0 & 0 & 0 & 0 & 0 & 0 & 0 & 0 & 0 & 0 & 0 \\
 1 & 0 & 0 & 0 & 0 & 0 & 0 & 0 & 0 & 0 & 0 \\
 0 & 0 & 0 & 0 & 0 & 0 & 0 & 0 & 0 & 0 & 0 \\
 0 & 0 & 0 & 0 & 0 & 0 & 0 & 0 & 0 & 0 & 0 \\
 0 & 0 & 0 & 0 & 0 & 0 & 0 & 0 & 0 & 0 & 0 \\
 0 & 0 & 0 & 0 & 0 & 0 & 0 & 0 & 0 & 0 & 0 \\
 0 & 0 & 0 & 0 & 0 & 0 & 0 & 0 & 0 & 0 & 0 \\
 0 & 0 & 0 & 0 & 0 & 0 & 0 & 0 & 0 & 0 & 0 \\
 0 & 0 & 0 & 0 & 0 & 0 & 0 & 0 & 0 & 0 & 0 \\
 0 & 0 & 0 & 0 & 0 & 0 & 0 & 0 & 0 & 0 & 0
\end{array}
\right)\ee 
\be \label{L_1R_3L_3R_1}
{\rm L}_1 *{\rm R}_4 + {\rm L}_4 *{\rm R}_1=\left(
\begin{array}{ccccccccccc}
 0 & 0 & 0 & -1 & 0 & 0 & 0 & 0 & 0 & 0 & 0 \\
 0 & 0 & 0 & 0 & 0 & 0 & 0 & 0 & 0 & 0 & 0 \\
 0 & 0 & 0 & 0 & 0 & 0 & 0 & 0 & 0 & 0 & 0 \\
 -1 & 0 & 0 & 0 & 0 & 0 & 0 & 0 & 0 & 0 & 0 \\
 0 & 0 & 0 & 0 & 0 & 0 & 0 & 0 & 0 & 0 & 0 \\
 0 & 0 & 0 & 0 & 0 & 0 & 0 & 0 & 0 & 0 & 0 \\
 0 & 0 & 0 & 0 & 0 & 0 & 0 & 0 & 0 & 0 & 0 \\
 0 & 0 & 0 & 0 & 0 & 0 & 0 & 0 & 0 & 0 & 0 \\
 0 & 0 & 0 & 0 & 0 & 0 & 0 & 0 & 0 & 0 & 0 \\
 0 & 0 & 0 & 0 & 0 & 0 & 0 & 0 & 0 & 0 & 0 \\
 0 & 0 & 0 & 0 & 0 & 0 & 0 & 0 & 0 & 0 & 0
\end{array}
\right)\ee 
\be {\rm L}_1 *{\rm R}_5 + {\rm L}_5 *{\rm R}_1=\left(
\begin{array}{ccccccccccc}
 0 & 0 & 0 & 0 & 0 & 0 & 0 & 0 & 0 & 0 & 0 \\
 0 & 0 & 0 & 0 & 0 & 0 & 0 & 1 & 0 & 0 & 0 \\
 0 & 0 & 0 & 0 & 0 & 0 & 0 & 0 & 0 & 0 & 0 \\
 0 & 0 & 0 & 0 & 0 & 0 & 0 & 0 & 0 & 1 & 0 \\
 0 & 0 & 0 & 0 & 0 & 0 & 0 & 0 & 0 & 0 & 0 \\
 0 & 0 & 0 & 0 & 0 & 0 & 0 & 0 & 0 & 0 & 0 \\
 0 & 0 & 0 & 0 & 0 & 0 & 0 & 0 & 0 & 0 & 0 \\
 0 & 1 & 0 & 0 & 0 & 0 & 0 & 0 & 0 & 0 & 0 \\
 0 & 0 & 0 & 0 & 0 & 0 & 0 & 0 & 0 & 0 & 0 \\
 0 & 0 & 0 & 1 & 0 & 0 & 0 & 0 & 0 & 0 & 0 \\
 0 & 0 & 0 & 0 & 0 & 0 & 0 & 0 & 0 & 0 & 0
\end{array}
\right)\ee 
\be\label{L_2R_2} {\rm L}_2*{\rm R}_2=\left(
\begin{array}{ccccccccccc}
 1 & 0 & 0 & 0 & 0 & 0 & 0 & 0 & 0 & 0 & 0 \\
 0 & 1 & 0 & 0 & 0 & 0 & 0 & 0 & 0 & 0 & 0 \\
 0 & 0 & 1 & 0 & 0 & 0 & 0 & 0 & 0 & 0 & 0 \\
 0 & 0 & 0 & 0 & 0 & 0 & 0 & 0 & 0 & 0 & 0 \\
 0 & 0 & 0 & 0 & 0 & 0 & 0 & 0 & 0 & 0 & 0 \\
 0 & 0 & 0 & 0 & 0 & 0 & 0 & 0 & 0 & 0 & 0 \\
 0 & 0 & 0 & 0 & 0 & 0 & 1 & 0 & 0 & 0 & 0 \\
 0 & 0 & 0 & 0 & 0 & 0 & 0 & 1 & 0 & 0 & 0 \\
 0 & 0 & 0 & 0 & 0 & 0 & 0 & 0 & 1 & 0 & 0 \\
 0 & 0 & 0 & 0 & 0 & 0 & 0 & 0 & 0 & 1 & 0 \\
 0 & 0 & 0 & 0 & 0 & 0 & 0 & 0 & 0 & 0 & 1
\end{array}
\right)\ee 
\be \label{L_2R_3L_3R_2}{\rm L}_2 *{\rm R}_3 + {\rm L}_3 *{\rm R}_2=\left(
\begin{array}{ccccccccccc}
 0 & 0 & 0 & 0 & 0 & 1 & 0 & 0 & 0 & 0 & 0 \\
 0 & 0 & 0 & 0 & 0 & 0 & 0 & 0 & 0 & 0 & 0 \\
 0 & 0 & 0 & 0 & 0 & 0 & 0 & 0 & 0 & 0 & 0 \\
 0 & 0 & 0 & 0 & 0 & 0 & -1 & 0 & 0 & 0 & 0 \\
 0 & 0 & 0 & 0 & 0 & 0 & 0 & 0 & 0 & 0 & 0 \\
 1 & 0 & 0 & 0 & 0 & 0 & 0 & 0 & 0 & 0 & 0 \\
 0 & 0 & 0 & -1 & 0 & 0 & 0 & 0 & 0 & 0 & 0 \\
 0 & 0 & 0 & 0 & 0 & 0 & 0 & 0 & 0 & 0 & 0 \\
 0 & 0 & 0 & 0 & 0 & 0 & 0 & 0 & 0 & 0 & 0 \\
 0 & 0 & 0 & 0 & 0 & 0 & 0 & 0 & 0 & 0 & 0 \\
 0 & 0 & 0 & 0 & 0 & 0 & 0 & 0 & 0 & 0 & 0
\end{array}
\right)\ee 
\be  \label{L_2R_4L_4R_2} {\rm L}_2*{\rm R}_4+{\rm L}_4*{\rm R}_2=\left(
\begin{array}{ccccccccccc}
 0 & 0 & 0 & 0 & 0 & 0 & 0 & 0 & 0 & 0 & 0 \\
 0 & 0 & 0 & 1 & 0 & 0 & 0 & 0 & 0 & 0 & 0 \\
 0 & 0 & 0 & 0 & 0 & 0 & 0 & 0 & 0 & 0 & 0 \\
 0 & 1 & 0 & 0 & 0 & 0 & 0 & 0 & 0 & 0 & 0 \\
 0 & 0 & 0 & 0 & 0 & 0 & 0 & 0 & 0 & 0 & 0 \\
 0 & 0 & 0 & 0 & 0 & 0 & 0 & 0 & 0 & 0 & 0 \\
 0 & 0 & 0 & 0 & 0 & 0 & 0 & 0 & 0 & 0 & 0 \\
 0 & 0 & 0 & 0 & 0 & 0 & 0 & 0 & 0 & 0 & 0 \\
 0 & 0 & 0 & 0 & 0 & 0 & 0 & 0 & 0 & 0 & 0 \\
 0 & 0 & 0 & 0 & 0 & 0 & 0 & 0 & 0 & 0 & 0 \\
 0 & 0 & 0 & 0 & 0 & 0 & 0 & 0 & 0 & 0 & 0
\end{array}
\right)\ee 
\be {\rm L}_2 *{\rm R}_5 + {\rm L}_5 *{\rm R}_2=\left(
\begin{array}{ccccccccccc}
 0 & 0 & 0 & 0 & 0 & 0 & 0 & 0 & 0 & 0 & 0 \\
 0 & 0 & 0 & 0 & -1 & 0 & 0 & 0 & 0 & 0 & 0 \\
 0 & 0 & 0 & 0 & 0 & 0 & 0 & 0 & 0 & 0 & 0 \\
 0 & 0 & 0 & 0 & 0 & 0 & 0 & 0 & 0 & 0 & 0 \\
 0 & -1 & 0 & 0 & 0 & 0 & 0 & 0 & 0 & 0 & 0 \\
 0 & 0 & 0 & 0 & 0 & 0 & 0 & 0 & 0 & 0 & 0 \\
 0 & 0 & 0 & 0 & 0 & 0 & 0 & 0 & 0 & 0 & 0 \\
 0 & 0 & 0 & 0 & 0 & 0 & 0 & 0 & 0 & 0 & 0 \\
 0 & 0 & 0 & 0 & 0 & 0 & 0 & 0 & 0 & 0 & 0 \\
 0 & 0 & 0 & 0 & 0 & 0 & 0 & 0 & 0 & 0 & 0 \\
 0 & 0 & 0 & 0 & 0 & 0 & 0 & 0 & 0 & 0 & 0
\end{array}
\right)\ee 
\be \label{L_3R_3} {\rm L}_3 *{\rm R}_3=\left(
\begin{array}{ccccccccccc}
 0 & 0 & 0 & 0 & 0 & 0 & 0 & 0 & 0 & 0 & 0 \\
 0 & 1 & 0 & 0 & 0 & 0 & 0 & 0 & 0 & 0 & 0 \\
 0 & 0 & 1 & 0 & 0 & 0 & 0 & 0 & 0 & 0 & 0 \\
 0 & 0 & 0 & 1 & 0 & 0 & 0 & 0 & 0 & 0 & 0 \\
 0 & 0 & 0 & 0 & 0 & 0 & 0 & 0 & 0 & 0 & 0 \\
 0 & 0 & 0 & 0 & 0 & 1 & 0 & 0 & 0 & 0 & 0 \\
 0 & 0 & 0 & 0 & 0 & 0 & 0 & 0 & 0 & 0 & 0 \\
 0 & 0 & 0 & 0 & 0 & 0 & 0 & 1 & 0 & 0 & 0 \\
 0 & 0 & 0 & 0 & 0 & 0 & 0 & 0 & 1 & 0 & 0 \\
 0 & 0 & 0 & 0 & 0 & 0 & 0 & 0 & 0 & 1 & 0 \\
 0 & 0 & 0 & 0 & 0 & 0 & 0 & 0 & 0 & 0 & 1
\end{array}
\right)\ee 
\be\label{L_3R_4L_4R_3} {\rm L}_3*{\rm R}_4+{\rm L}_4*{\rm R}_3= \left(
\begin{array}{ccccccccccc}
 0 & 0 & 0 & 0 & 0 & 0 & 0 & 0 & 0 & 0 & 0 \\
 0 & 0 & 0 & 0 & 0 & 0 & 1 & 0 & 0 & 0 & 0 \\
 0 & 0 & 0 & 0 & 0 & 0 & 0 & 0 & 0 & 0 & 0 \\
 0 & 0 & 0 & 0 & 0 & 0 & 0 & 0 & 0 & 0 & 0 \\
 0 & 0 & 0 & 0 & 0 & 0 & 0 & 0 & -1 & 0 & 0 \\
 0 & 0 & 0 & 0 & 0 & 0 & 0 & 0 & 0 & 0 & 0 \\
 0 & 1 & 0 & 0 & 0 & 0 & 0 & 0 & 0 & 0 & 0 \\
 0 & 0 & 0 & 0 & 0 & 0 & 0 & 0 & 0 & 0 & 0 \\
 0 & 0 & 0 & 0 & -1 & 0 & 0 & 0 & 0 & 0 & 0 \\
 0 & 0 & 0 & 0 & 0 & 0 & 0 & 0 & 0 & 0 & 0 \\
 0 & 0 & 0 & 0 & 0 & 0 & 0 & 0 & 0 & 0 & 0
\end{array}
\right)\ee 
\be {\rm L}_3 *{\rm R}_5 + {\rm L}_5 *{\rm R}_3=\left(
\begin{array}{ccccccccccc}
 0 & 0 & 0 & 0 & 0 & 0 & 0 & 0 & 0 & 0 & 0 \\
 0 & 0 & 0 & 0 & 0 & 0 & 0 & 0 & 0 & 0 & 0 \\
 0 & 0 & 0 & 0 & 1 & 0 & 0 & 0 & 0 & 0 & 0 \\
 0 & 0 & 0 & 0 & 0 & 0 & 0 & 0 & 0 & 0 & 0 \\
 0 & 0 & 1 & 0 & 0 & 0 & 0 & 0 & 0 & 0 & 0 \\
 0 & 0 & 0 & 0 & 0 & 0 & 0 & 0 & 0 & 0 & 0 \\
 0 & 0 & 0 & 0 & 0 & 0 & 0 & 0 & 0 & 0 & 0 \\
 0 & 0 & 0 & 0 & 0 & 0 & 0 & 0 & 0 & 0 & 0 \\
 0 & 0 & 0 & 0 & 0 & 0 & 0 & 0 & 0 & 0 & 0 \\
 0 & 0 & 0 & 0 & 0 & 0 & 0 & 0 & 0 & 0 & 0 \\
 0 & 0 & 0 & 0 & 0 & 0 & 0 & 0 & 0 & 0 & 0
\end{array}
\right)\ee 
\be\label{L_4R_4} {\rm L}_4*{\rm R}_4 = \left(
\begin{array}{ccccccccccc}
 0 & 0 & 0 & 0 & 0 & 0 & 0 & 0 & 0 & 0 & 0 \\
 0 & 0 & 0 & 0 & 0 & 0 & 0 & 0 & 0 & 0 & 0 \\
 0 & 0 & 1 & 0 & 0 & 0 & 0 & 0 & 0 & 0 & 0 \\
 0 & 0 & 0 & 1 & 0 & 0 & 0 & 0 & 0 & 0 & 0 \\
 0 & 0 & 0 & 0 & 1 & 0 & 0 & 0 & 0 & 0 & 0 \\
 0 & 0 & 0 & 0 & 0 & 1 & 0 & 0 & 0 & 0 & 0 \\
 0 & 0 & 0 & 0 & 0 & 0 & 1 & 0 & 0 & 0 & 0 \\
 0 & 0 & 0 & 0 & 0 & 0 & 0 & 1 & 0 & 0 & 0 \\
 0 & 0 & 0 & 0 & 0 & 0 & 0 & 0 & 0 & 0 & 0 \\
 0 & 0 & 0 & 0 & 0 & 0 & 0 & 0 & 0 & 1 & 0 \\
 0 & 0 & 0 & 0 & 0 & 0 & 0 & 0 & 0 & 0 & 1
\end{array}
\right)\ee 
\be {\rm L}_4 *{\rm R}_5 + {\rm L}_5 *{\rm R}_4=\left(
\begin{array}{ccccccccccc}
 0 & 0 & 0 & 0 & 0 & 0 & 0 & 0 & 0 & 1 & 0 \\
 0 & 0 & 0 & 0 & 0 & 0 & 0 & 0 & 0 & 0 & 0 \\
 0 & 0 & 0 & 0 & 0 & 0 & 0 & 0 & 1 & 0 & 0 \\
 0 & 0 & 0 & 0 & 0 & 0 & 0 & 0 & 0 & 0 & 0 \\
 0 & 0 & 0 & 0 & 0 & 0 & 0 & 0 & 0 & 0 & 0 \\
 0 & 0 & 0 & 0 & 0 & 0 & 0 & 0 & 0 & 0 & 0 \\
 0 & 0 & 0 & 0 & 0 & 0 & 0 & 0 & 0 & 0 & 0 \\
 0 & 0 & 0 & 0 & 0 & 0 & 0 & 0 & 0 & 0 & 0 \\
 0 & 0 & 1 & 0 & 0 & 0 & 0 & 0 & 0 & 0 & 0 \\
 1 & 0 & 0 & 0 & 0 & 0 & 0 & 0 & 0 & 0 & 0 \\
 0 & 0 & 0 & 0 & 0 & 0 & 0 & 0 & 0 & 0 & 0
\end{array}
\right)\ee 
\be\label{L_5R_5} {\rm L}_5*{\rm R}_5 = \left(
\begin{array}{ccccccccccc}
 1 & 0 & 0 & 0 & 0 & 0 & 0 & 0 & 0 & 0 & 0 \\
 0 & 0 & 0 & 0 & 0 & 0 & 0 & 0 & 0 & 0 & 0 \\
 0 & 0 & 0 & 0 & 0 & 0 & 0 & 0 & 0 & 0 & 0 \\
 0 & 0 & 0 & 1 & 0 & 0 & 0 & 0 & 0 & 0 & 0 \\
 0 & 0 & 0 & 0 & 1 & 0 & 0 & 0 & 0 & 0 & 0 \\
 0 & 0 & 0 & 0 & 0 & 1 & 0 & 0 & 0 & 0 & 0 \\
 0 & 0 & 0 & 0 & 0 & 0 & 1 & 0 & 0 & 0 & 0 \\
 0 & 0 & 0 & 0 & 0 & 0 & 0 & 1 & 0 & 0 & 0 \\
 0 & 0 & 0 & 0 & 0 & 0 & 0 & 0 & 1 & 0 & 0 \\
 0 & 0 & 0 & 0 & 0 & 0 & 0 & 0 & 0 & 0 & 0 \\
 0 & 0 & 0 & 0 & 0 & 0 & 0 & 0 & 0 & 0 & 1
\end{array}
\right)\ee \\
\\

\be \label{L_1R_1} {\rm R}_1*{\rm L}_1=\left(
\begin{array}{ccccccccccc}
 1 & 0 & 0 & 0 & 0 & 0 & 0 & 0 & 0 & 0 & 0 \\
 0 & 1 & 0 & 0 & 0 & 0 & 0 & 0 & 0 & 0 & 0 \\
 0 & 0 & 1 & 0 & 0 & 0 & 0 & 0 & 0 & 0 & 0 \\
 0 & 0 & 0 & 1 & 0 & 0 & 0 & 0 & 0 & 0 & 0 \\
 0 & 0 & 0 & 0 & 0 & 0 & 0 & 0 & 0 & 0 & 0 \\
 0 & 0 & 0 & 0 & 0 & 1 & 0 & 0 & 0 & 0 & 0 \\
 0 & 0 & 0 & 0 & 0 & 0 & 1 & 0 & 0 & 0 & 0 \\
 0 & 0 & 0 & 0 & 0 & 0 & 0 & 1 & 0 & 0 & 0 \\
 0 & 0 & 0 & 0 & 0 & 0 & 0 & 0 & 0 & 0 & 0 \\
 0 & 0 & 0 & 0 & 0 & 0 & 0 & 0 & 0 & 1 & 0 \\
 0 & 0 & 0 & 0 & 0 & 0 & 0 & 0 & 0 & 0 & 0
\end{array}
\right)\ee 
\be\label{R_1L_2R_2L_1} {\rm R}_1*{\rm L}_2 + {\rm R}_2*{\rm L}_1=\left(
\begin{array}{ccccccccccc}
 0 & 0 & 0 & 0 & 0 & 0 & 0 & 0 & 0 & 0 & 0 \\
 0 & 0 & 0 & 0 & 0 & 0 & 0 & 0 & 0 & 0 & 0 \\
 0 & 0 & 0 & 0 & 0 & 0 & 0 & 0 & 0 & 0 & 0 \\
 0 & 0 & 0 & 0 & 0 & 0 & 0 & 0 & 0 & 0 & 0 \\
 0 & 0 & 0 & 0 & 0 & 0 & -1 & 0 & 0 & 0 & 0 \\
 0 & 0 & 0 & 0 & 0 & 0 & 0 & 0 & 0 & 0 & 1 \\
 0 & 0 & 0 & 0 & -1 & 0 & 0 & 0 & 0 & 0 & 0 \\
 0 & 0 & 0 & 0 & 0 & 0 & 0 & 0 & 0 & 0 & 0 \\
 0 & 0 & 0 & 0 & 0 & 0 & 0 & 0 & 0 & 0 & 0 \\
 0 & 0 & 0 & 0 & 0 & 0 & 0 & 0 & 0 & 0 & 0 \\
 0 & 0 & 0 & 0 & 0 & 1 & 0 & 0 & 0 & 0 & 0
\end{array}
\right)\ee 
\be\label{R_1L_3R_3L_1} {\rm R}_1*{\rm L}_3+{\rm R}_3 *{\rm L}_1= \left(
\begin{array}{ccccccccccc}
 0 & 0 & 0 & 0 & 0 & 0 & 0 & 0 & 0 & 0 & 0 \\
 0 & 0 & 0 & 0 & 0 & 0 & 0 & 0 & 0 & 0 & 0 \\
 0 & 0 & 0 & 0 & 0 & 0 & 0 & 0 & 0 & 0 & 0 \\
 0 & 0 & 0 & 0 & 0 & 0 & 0 & 0 & 0 & 0 & 0 \\
 0 & 0 & 0 & 0 & 0 & 0 & 0 & 0 & 0 & 0 & 0 \\
 0 & 0 & 0 & 0 & 0 & 0 & 0 & 0 & 0 & 0 & 0 \\
 0 & 0 & 0 & 0 & 0 & 0 & 0 & 0 & 0 & 0 & 0 \\
 0 & 0 & 0 & 0 & 0 & 0 & 0 & 0 & 0 & 0 & 1 \\
 0 & 0 & 0 & 0 & 0 & 0 & 0 & 0 & 0 & 0 & 0 \\
 0 & 0 & 0 & 0 & 0 & 0 & 0 & 0 & 0 & 0 & 0 \\
 0 & 0 & 0 & 0 & 0 & 0 & 0 & 1 & 0 & 0 & 0
\end{array}
\right)\ee 
\be  \label{R_1L_4R_4L_1}{\rm R}_1 *{\rm L}_4 + {\rm R}_4 *{\rm L}_1=\left(
\begin{array}{ccccccccccc}
 0 & 0 & 0 & 0 & 0 & 0 & 0 & 0 & 0 & 0 & 0 \\
 0 & 0 & 0 & 0 & 0 & 0 & 0 & 0 & 0 & 0 & 0 \\
 0 & 0 & 0 & 0 & 0 & 0 & 0 & 0 & 0 & 0 & 0 \\
 0 & 0 & 0 & 0 & 0 & 0 & 0 & 0 & 0 & 0 & 0 \\
 0 & 0 & 0 & 0 & 0 & 0 & 0 & 0 & 0 & 0 & 0 \\
 0 & 0 & 0 & 0 & 0 & 0 & 0 & 0 & 0 & 0 & 0 \\
 0 & 0 & 0 & 0 & 0 & 0 & 0 & 0 & 0 & 0 & 0 \\
 0 & 0 & 0 & 0 & 0 & 0 & 0 & 0 & -1 & 0 & 0 \\
 0 & 0 & 0 & 0 & 0 & 0 & 0 & -1 & 0 & 0 & 0 \\
 0 & 0 & 0 & 0 & 0 & 0 & 0 & 0 & 0 & 0 & 0 \\
 0 & 0 & 0 & 0 & 0 & 0 & 0 & 0 & 0 & 0 & 0
\end{array}
\right)\ee 
\be {\rm R}_1*{\rm L}_5 + {\rm R}_5 *{\rm L}_1=\left(
\begin{array}{ccccccccccc}
 0 & 0 & 0 & 0 & 0 & 0 & 0 & 0 & 0 & 0 & 0 \\
 0 & 0 & 0 & 0 & 0 & 0 & 0 & 0 & 0 & 0 & 0 \\
 0 & 0 & 0 & 0 & 0 & 0 & 0 & 0 & 1 & 0 & 0 \\
 0 & 0 & 0 & 0 & 0 & 0 & 0 & 0 & 0 & 0 & 0 \\
 0 & 0 & 0 & 0 & 0 & 0 & 0 & 0 & 0 & 1 & 0 \\
 0 & 0 & 0 & 0 & 0 & 0 & 0 & 0 & 0 & 0 & 0 \\
 0 & 0 & 0 & 0 & 0 & 0 & 0 & 0 & 0 & 0 & 0 \\
 0 & 0 & 0 & 0 & 0 & 0 & 0 & 0 & 0 & 0 & 0 \\
 0 & 0 & 1 & 0 & 0 & 0 & 0 & 0 & 0 & 0 & 0 \\
 0 & 0 & 0 & 0 & 1 & 0 & 0 & 0 & 0 & 0 & 0 \\
 0 & 0 & 0 & 0 & 0 & 0 & 0 & 0 & 0 & 0 & 0
\end{array}
\right)\ee 
\be\label{R_2L_2} {\rm R}_2*{\rm L}_2=\left(
\begin{array}{ccccccccccc}
 1 & 0 & 0 & 0 & 0 & 0 & 0 & 0 & 0 & 0 & 0 \\
 0 & 1 & 0 & 0 & 0 & 0 & 0 & 0 & 0 & 0 & 0 \\
 0 & 0 & 1 & 0 & 0 & 0 & 0 & 0 & 0 & 0 & 0 \\
 0 & 0 & 0 & 1 & 0 & 0 & 0 & 0 & 0 & 0 & 0 \\
 0 & 0 & 0 & 0 & 1 & 0 & 0 & 0 & 0 & 0 & 0 \\
 0 & 0 & 0 & 0 & 0 & 0 & 0 & 0 & 0 & 0 & 0 \\
 0 & 0 & 0 & 0 & 0 & 0 & 0 & 0 & 0 & 0 & 0 \\
 0 & 0 & 0 & 0 & 0 & 0 & 0 & 1 & 0 & 0 & 0 \\
 0 & 0 & 0 & 0 & 0 & 0 & 0 & 0 & 0 & 0 & 0 \\
 0 & 0 & 0 & 0 & 0 & 0 & 0 & 0 & 0 & 1 & 0 \\
 0 & 0 & 0 & 0 & 0 & 0 & 0 & 0 & 0 & 0 & 1
\end{array}
\right)\ee 
\be\label{R_2L_3R_3L_2} {\rm R}_2*{\rm L}_3 + {\rm R}_3*{\rm L}_2=\left(
\begin{array}{ccccccccccc}
 0 & 0 & 0 & 0 & 0 & 0 & 0 & 0 & 0 & 0 & 0 \\
 0 & 0 & 0 & 0 & 0 & 0 & 0 & 0 & 0 & 0 & 0 \\
 0 & 0 & 0 & 0 & 0 & 0 & 0 & 0 & 0 & 0 & 0 \\
 0 & 0 & 0 & 0 & 0 & 0 & 0 & 0 & -1 & 0 & 0 \\
 0 & 0 & 0 & 0 & 0 & 0 & 0 & 0 & 0 & 0 & 0 \\
 0 & 0 & 0 & 0 & 0 & 0 & 0 & -1 & 0 & 0 & 0 \\
 0 & 0 & 0 & 0 & 0 & 0 & 0 & 0 & 0 & 0 & 0 \\
 0 & 0 & 0 & 0 & 0 & -1 & 0 & 0 & 0 & 0 & 0 \\
 0 & 0 & 0 & -1 & 0 & 0 & 0 & 0 & 0 & 0 & 0 \\
 0 & 0 & 0 & 0 & 0 & 0 & 0 & 0 & 0 & 0 & 0 \\
 0 & 0 & 0 & 0 & 0 & 0 & 0 & 0 & 0 & 0 & 0
\end{array}
\right)\ee 
\be \label{R_2L_4R_4L_2} {\rm R}_2*{\rm L}_4 + {\rm R}_4*{\rm L}_2=\left(
\begin{array}{ccccccccccc}
 0 & 0 & 0 & 0 & 0 & 0 & 0 & 0 & 0 & 0 & 0 \\
 0 & 0 & 0 & 0 & 0 & 0 & 0 & 0 & 0 & 0 & 0 \\
 0 & 0 & 0 & 0 & 0 & 0 & 0 & 0 & 0 & 0 & 0 \\
 0 & 0 & 0 & 0 & 0 & 0 & 0 & 0 & 0 & 0 & 0 \\
 0 & 0 & 0 & 0 & 0 & 0 & 0 & 0 & 0 & 0 & 0 \\
 0 & 0 & 0 & 0 & 0 & 0 & 0 & 0 & 0 & 0 & 0 \\
 0 & 0 & 0 & 0 & 0 & 0 & 0 & 0 & 0 & 0 & 0 \\
 0 & 0 & 0 & 0 & 0 & 0 & 0 & 0 & 0 & 0 & 0 \\
 0 & 0 & 0 & 0 & 0 & 0 & 0 & 0 & 0 & 1 & 0 \\
 0 & 0 & 0 & 0 & 0 & 0 & 0 & 0 & 1 & 0 & 0 \\
 0 & 0 & 0 & 0 & 0 & 0 & 0 & 0 & 0 & 0 & 0
\end{array}
\right)\ee 
\be {\rm R}_2*{\rm L}_5 + {\rm R}_5 *{\rm L}_2=\left(
\begin{array}{ccccccccccc}
 0 & 0 & 0 & 0 & 0 & 0 & 0 & 0 & 0 & 0 & 0 \\
 0 & 0 & 0 & 0 & 0 & 0 & 0 & 0 & 0 & 0 & 0 \\
 0 & 0 & 0 & 0 & 0 & 0 & 0 & 0 & 0 & 0 & 0 \\
 0 & 0 & 0 & 0 & 0 & 0 & 0 & 0 & 0 & 0 & 0 \\
 0 & 0 & 0 & 0 & 0 & 0 & 0 & 0 & 0 & 0 & 0 \\
 0 & 0 & 0 & 0 & 0 & 0 & 0 & 0 & 0 & 0 & 0 \\
 0 & 0 & 0 & 0 & 0 & 0 & 0 & 0 & 0 & 1 & 0 \\
 0 & 0 & 0 & 0 & 0 & 0 & 0 & 0 & 0 & 0 & 0 \\
 0 & 0 & 0 & 0 & 0 & 0 & 0 & 0 & 0 & 0 & 0 \\
 0 & 0 & 0 & 0 & 0 & 0 & 1 & 0 & 0 & 0 & 0 \\
 0 & 0 & 0 & 0 & 0 & 0 & 0 & 0 & 0 & 0 & 0
\end{array}
\right)\ee 
\be \label{R_3L_3}{\rm R}_3*{\rm L}_3=\left(
\begin{array}{ccccccccccc}
 1 & 0 & 0 & 0 & 0 & 0 & 0 & 0 & 0 & 0 & 0 \\
 0 & 1 & 0 & 0 & 0 & 0 & 0 & 0 & 0 & 0 & 0 \\
 0 & 0 & 1 & 0 & 0 & 0 & 0 & 0 & 0 & 0 & 0 \\
 0 & 0 & 0 & 0 & 0 & 0 & 0 & 0 & 0 & 0 & 0 \\
 0 & 0 & 0 & 0 & 1 & 0 & 0 & 0 & 0 & 0 & 0 \\
 0 & 0 & 0 & 0 & 0 & 1 & 0 & 0 & 0 & 0 & 0 \\
 0 & 0 & 0 & 0 & 0 & 0 & 0 & 0 & 0 & 0 & 0 \\
 0 & 0 & 0 & 0 & 0 & 0 & 0 & 0 & 0 & 0 & 0 \\
 0 & 0 & 0 & 0 & 0 & 0 & 0 & 0 & 1 & 0 & 0 \\
 0 & 0 & 0 & 0 & 0 & 0 & 0 & 0 & 0 & 1 & 0 \\
 0 & 0 & 0 & 0 & 0 & 0 & 0 & 0 & 0 & 0 & 1
\end{array}
\right)\ee 
\be  \label{R_3L_4R_4L_3}{\rm R}_3*{\rm L}_4+{\rm R}_4*{\rm L}_3=\left(
\begin{array}{ccccccccccc}
 0 & 0 & 0 & 0 & 0 & 0 & 0 & 0 & 0 & 0 & 0 \\
 0 & 0 & 0 & 0 & 0 & 0 & 1 & 0 & 0 & 0 & 0 \\
 0 & 0 & 0 & 0 & 0 & 0 & 0 & 0 & 0 & 0 & 0 \\
 0 & 0 & 0 & 0 & 0 & 0 & 0 & 0 & 0 & 1 & 0 \\
 0 & 0 & 0 & 0 & 0 & 0 & 0 & 0 & 0 & 0 & 0 \\
 0 & 0 & 0 & 0 & 0 & 0 & 0 & 0 & 0 & 0 & 0 \\
 0 & 1 & 0 & 0 & 0 & 0 & 0 & 0 & 0 & 0 & 0 \\
 0 & 0 & 0 & 0 & 0 & 0 & 0 & 0 & 0 & 0 & 0 \\
 0 & 0 & 0 & 0 & 0 & 0 & 0 & 0 & 0 & 0 & 0 \\
 0 & 0 & 0 & 1 & 0 & 0 & 0 & 0 & 0 & 0 & 0 \\
 0 & 0 & 0 & 0 & 0 & 0 & 0 & 0 & 0 & 0 & 0
\end{array}
\right)\ee 
\be {\rm R}_3*{\rm L}_5 + {\rm R}_5 *{\rm L}_3=\left(
\begin{array}{ccccccccccc}
 0 & 0 & 0 & 0 & 0 & 0 & 0 & 0 & 0 & 0 & 0 \\
 0 & 0 & 0 & 0 & 0 & 0 & 0 & 0 & 0 & 0 & 0 \\
 0 & 0 & 0 & 0 & 0 & 0 & 0 & 0 & 0 & 0 & 0 \\
 0 & 0 & 0 & 0 & 0 & 0 & 0 & 0 & 0 & 0 & 0 \\
 0 & 0 & 0 & 0 & 0 & 0 & 0 & 0 & 0 & 0 & 0 \\
 0 & 0 & 0 & 0 & 0 & 0 & 0 & 0 & 0 & 0 & 0 \\
 0 & 0 & 0 & 0 & 0 & 0 & 0 & 0 & 0 & 0 & -1 \\
 0 & 0 & 0 & 0 & 0 & 0 & 0 & 0 & 0 & 0 & 0 \\
 0 & 0 & 0 & 0 & 0 & 0 & 0 & 0 & 0 & 0 & 0 \\
 0 & 0 & 0 & 0 & 0 & 0 & 0 & 0 & 0 & 0 & 0 \\
 0 & 0 & 0 & 0 & 0 & 0 & -1 & 0 & 0 & 0 & 0
\end{array}
\right)\ee 
\be  \label{R4L4}{\rm R}_4*{\rm L}_4=\left(
\begin{array}{ccccccccccc}
 1 & 0 & 0 & 0 & 0 & 0 & 0 & 0 & 0 & 0 & 0 \\
 0 & 0 & 0 & 0 & 0 & 0 & 0 & 0 & 0 & 0 & 0 \\
 0 & 0 & 1 & 0 & 0 & 0 & 0 & 0 & 0 & 0 & 0 \\
 0 & 0 & 0 & 1 & 0 & 0 & 0 & 0 & 0 & 0 & 0 \\
 0 & 0 & 0 & 0 & 1 & 0 & 0 & 0 & 0 & 0 & 0 \\
 0 & 0 & 0 & 0 & 0 & 1 & 0 & 0 & 0 & 0 & 0 \\
 0 & 0 & 0 & 0 & 0 & 0 & 1 & 0 & 0 & 0 & 0 \\
 0 & 0 & 0 & 0 & 0 & 0 & 0 & 0 & 0 & 0 & 0 \\
 0 & 0 & 0 & 0 & 0 & 0 & 0 & 0 & 1 & 0 & 0 \\
 0 & 0 & 0 & 0 & 0 & 0 & 0 & 0 & 0 & 0 & 0 \\
 0 & 0 & 0 & 0 & 0 & 0 & 0 & 0 & 0 & 0 & 1
\end{array}
\right)\ee 
\be {\rm R}_4*{\rm L}_5 + {\rm R}_5 *{\rm L}_4=\left(
\begin{array}{ccccccccccc}
 0 & 0 & 0 & 0 & 0 & 0 & 0 & 0 & 0 & 0 & 0 \\
 0 & 0 & 0 & 0 & 0 & 0 & 0 & 0 & 0 & 0 & 1 \\
 0 & 0 & 0 & 0 & 0 & 0 & 0 & 1 & 0 & 0 & 0 \\
 0 & 0 & 0 & 0 & 0 & 0 & 0 & 0 & 0 & 0 & 0 \\
 0 & 0 & 0 & 0 & 0 & 0 & 0 & 0 & 0 & 0 & 0 \\
 0 & 0 & 0 & 0 & 0 & 0 & 0 & 0 & 0 & 0 & 0 \\
 0 & 0 & 0 & 0 & 0 & 0 & 0 & 0 & 0 & 0 & 0 \\
 0 & 0 & 1 & 0 & 0 & 0 & 0 & 0 & 0 & 0 & 0 \\
 0 & 0 & 0 & 0 & 0 & 0 & 0 & 0 & 0 & 0 & 0 \\
 0 & 0 & 0 & 0 & 0 & 0 & 0 & 0 & 0 & 0 & 0 \\
 0 & 1 & 0 & 0 & 0 & 0 & 0 & 0 & 0 & 0 & 0
\end{array}
\right)\ee 
\be  \label{R5L5}{\rm R}_5*{\rm L}_5=\left(
\begin{array}{ccccccccccc}
 1 & 0 & 0 & 0 & 0 & 0 & 0 & 0 & 0 & 0 & 0 \\
 0 & 1 & 0 & 0 & 0 & 0 & 0 & 0 & 0 & 0 & 0 \\
 0 & 0 & 0 & 0 & 0 & 0 & 0 & 0 & 0 & 0 & 0 \\
 0 & 0 & 0 & 1 & 0 & 0 & 0 & 0 & 0 & 0 & 0 \\
 0 & 0 & 0 & 0 & 1 & 0 & 0 & 0 & 0 & 0 & 0 \\
 0 & 0 & 0 & 0 & 0 & 1 & 0 & 0 & 0 & 0 & 0 \\
 0 & 0 & 0 & 0 & 0 & 0 & 1 & 0 & 0 & 0 & 0 \\
 0 & 0 & 0 & 0 & 0 & 0 & 0 & 1 & 0 & 0 & 0 \\
 0 & 0 & 0 & 0 & 0 & 0 & 0 & 0 & 1 & 0 & 0 \\
 0 & 0 & 0 & 0 & 0 & 0 & 0 & 0 & 0 & 0 & 0 \\
 0 & 0 & 0 & 0 & 0 & 0 & 0 & 0 & 0 & 0 & 0
\end{array}
\right)\ee 



\begin{thebibliography}{99}
\bibitem{Adinkra}
S.\ J.\ Gates, Jr.,   Michael Faux, 
``Adinkras: A Graphical Technology for Supersymmetric Representation Theory,"
Phys.\ Rev.\ {\bf {D71}} 065002 (2005), hep-th/0408004.

\bibitem{Cubes}
C.\ F.\ Doran, M.\ G.\ Faux, S.\ J.\ Gates, Jr., T.\ H\" ubsch, K.\ M.\ Iga, G.\ D.\ Landweber, and
R.\ L\. Miller, ``Topology Types of Adinkras and the Corresponding Representations of N-Extended Supersymmetry," UMD preprint \# 08-010, SUNY-O/667, arXiv:0806.0050 [hep-th]. 
 
 \bibitem{N2}
C.\ F.\ Doran, M.\ G.\ Faux, S.\ J.\ Gates, Jr., T.\ H\" ubsch, K.\ M.\ Iga, and G.\ D.\ Landweber,
``On the Matter of $\cal N$ $=$ 2 Matter,'' Phys.\ Lett.\ {\bf {B659}} 441 (2008), 
arXiv:0710.5245 [hep-th]; idem.\, ``Frames for supersymmetry,'' Int.\ J.\  Mod.\ Phys.\
{\bf {A24}} 065402 (2009), arXiv:0809.5279 [hep-th].
 
\bibitem{Genomics}
S.\ J.\ Gates, Jr., J.\ Gonzales, B.\ MacGregor, J.\ Parker, R.\ Polo-Sherk, V.\ G.\ J.\ Rodgers,
L.\ Wassink, ``4D, N = 1 Supersymmetry Genomics (I),''
arXiv:0902.3830 [hep-th], JHEP 0912:008, 2009.

\bibitem{ERIK} The suggestion that zonohedra be considered as the starting point for the possible construction of adinkras was made to SJG by Erik Dumaine .

\end{thebibliography}
\end{document}